\documentclass[12pt]{amsart}
\input{epsf}
\usepackage{amssymb,latexsym,cite}
\topskip  \numberwithin{equation}{section}
\newtheorem{theorem}{Theorem}[section]
\newtheorem{proposition}[theorem]{Proposition}

\newtheorem{remark}[theorem]{Remark}

\begin{document}

\pagenumbering{arabic}
\pagestyle{headings}

\newcommand{\DPB}[4]{P\beta_{#1}^{(#2)}(#3,#4)}

\title[Central factorial numbers associated with sequences of polynomials]{Central factorial numbers associated with sequences of polynomials}
\author{Dae San Kim}
\address{Department of Mathematics, Sogang University, Seoul 121-742, Republic of Korea}
\email{dskim@sogang.ac.kr}

\author{Taekyun Kim}
\address{Department of Mathematics, Kwangwoon University, Seoul 139-701, Republic of Korea}
\email{tkkim@kw.ac.kr}


\begin{abstract}
Many important special numbers appear in the expansions of some polynomials in terms of central factorials and vice versa, for example central factorial numbers, degenerate central factorial numbers, and central Lah numbers which are recently introduced. Here we generalize this to any sequence of polynomials. Let $\bold{P}=\left\{p_{n}(x)\right\}_{n=0}^{\infty}$ be a sequence of polynomials such that deg\,$p_{n}(x)=n, p_{0}(x)=1$. The aim of this paper is to study the central factorial numbers of the second associated with $\bold{P}$ and of the first kind associated with $\bold{P}$, in a unified and systematic way with the help of umbral calculus technique. The central factorial numbers associated with $\bold{P}$ enjoy orthogonality and inverse relations. We illustrate our results with many examples and obtain interesting orthogonality and inverse relations by applying such relations for the central factorial numbers associated with $\bold{P}$  to each of our examples.
\end{abstract}

\maketitle

\markboth{\centerline{\scriptsize Central factorial numbers associated with sequences of polynomials}}{\centerline{\scriptsize Dae San Kim, Taekyun Kim}}


\section{Introduction}
The central factorial numbers of the second $T_{2}(n,k)$ appear as the coefficients in the expansion of the powers of $x$ in terms of the central factorials (see \eqref{23B}, \eqref{2C}), whereas the central factorial numbers of the first kind $T_{1}(n,k)$ appear as the coefficients in the expansion of the central factorials in terms of powers of $x$ (see \eqref{1C}). Combinatorial interpretations of central factorial numbers of the second kind are as follows. $T_{2}(n,k)$ is equal to the number of partitions of the set $\left\{1,1',2,2',...,n,n'\right\}$ into $k$ disjoint nonempty subsets $A_{1},...,A_{k}$ such that, for each $1 \le j \le k$, if $i$ is the least integer such that either $i$ or $i'$ belongs to $V_{j}$ then $\left\{i,i'\right\}$ is a subset of $A_{j}$. For example, $T_{2}(3,2) = 5$. Indeed, they are given by $\left\{1,1',2,2'\right\}\left\{3,3'\right\},\left\{1,1',3,3'\right\}\left\{2,2'\right\}, \left\{1,1'\right\}\left\{2,2',3,3'\right\}$, $\left\{1,1',3\right\}\left\{2,2',3'\right\}$, and $\left\{1,1',3'\right\}\left\{2,2',3\right\}$. It is also known that $T_{2}(2n, 2n–- 2k)$ enumerates the number of ways to pace $k$ rooks on a 3$D$-triangle board of size $(n-1)$ (see [28]). Further, the number of ways to place $k$ rooks on a size $n$ triangle board is $T_{2}(n+1,n+1-k)$ (see [28]). Even though the central factorial numbers are less known than Stirling numbers, they are as important as Stirling numbers.\\
Let $\bold{P}=\left\{p_{n}(x)\right\}_{n=0}^{\infty}$ be a sequence of polynomials such that deg\,$p_{n}(x)=n, p_{0}(x)=1$. 
The aim of this paper is to study the central factorial numbers of the second $T_{2}(n,k;\bold{P})$ associated with $\bold{P}$, which are defined as the coefficients in the expansion of $p_{n}(x)$ in terms of the central factorials $x^{[k]}$. This is done with the help of umbral calculus technique. We illustrate our results with twenty examples. Moreover, we also introduce the central factorial numbers of the first kind $T_{1}(n,k;\bold{P})$ associated with $\bold{P}$, which are defined as the coefficients in the expansion of $x^{[n]}$ in terms of the polynomials $p_{k}(x)$.  Again, we study them by means of umbral calculus and illustrate our results with the same twenty examples for the central factorial numbers of the second kind.
We note here that the orthogonality and inverse relations for the ordinary central factorial numbers also hold for the central factorial numbers associated with sequences of polynomials $\bold{P}$. This observation applied to each individual case yields interesting orthogonality and inverse relations for the twenty examples, which are explicitly stated.  For example,  in the case of the sequence $\bold{P}=\left\{\langle{x \rangle}_{n}\right\}$ of rising factorials, we obtain:
\begin{align*}
&\sum_{k=l}^{n}\sum_{m=k}^{n}\sum_{j=l}^{k}T_{1}(n,m)(-1)^{m-k}S_{2}(m,k)(-1)^{k-j}S_{1}(k,j)T_{2}(j,l)=\delta_{n,l},\\
&\sum_{k=l}^{n}\sum_{j=k}^{n}\sum_{m=l}^{k}(-1)^{n-j}S_{1}(n,j)T_{2}(j,k)T_{1}(k,m)(-1)^{m-l}S_{2}(m,l)=\delta_{n,l}, \\
&a_{n}=\sum_{k=0}^{n}\sum_{l=k}^{n}(-1)^{n-l}S_{1}(n,l)T_{2}(l,k)c_{k} \Longleftrightarrow c_{n}=\sum_{k=0}^{n}\sum_{l=k}^{n}(-1)^{l-k}T_{1}(n,l)S_{2}(l,k)a_{k},\\
&a_{n}=\sum_{k=n}^{m}\sum_{l=n}^{k}(-1)^{k-l}S_{1}(k,l)T_{2}(l,n)c_{k} \Longleftrightarrow c_{n}=\sum_{k=n}^{m}\sum_{l=n}^{k}(-1)^{l-n}T_{1}(k,l)S_{2}(l,n)a_{k}, 
\end{align*}
where $S_{1}(n,k), S_{2}(n,k)$ are the Stirling numbers of the first kind and the Stirling numbers of the second kind, respectively.
As another example, in the case of the sequence $\bold{P}=\left\{L_{n}(x)\right\}$ of Laguerre polynomials of order -1,  we have:
\begin{align*}
&\sum_{k=l}^{n}\sum_{m=k}^{n}\sum_{j=l}^{k}(-1)^{k-j}T_{1}(n,m)L(m,k)L(k,j)T_{2}(j,l)=\delta_{n,l},\\ 
&\sum_{k=l}^{n}\sum_{j=k}^{n}\sum_{m=l}^{k}(-1)^{j-l}L(n,j)T_{2}(j,k)T_{1}(k,m)L(m,l)=\delta_{n,l}, \\
&a_{n}=\sum_{k=0}^{n}\sum_{l=k}^{n}(-1)^{l}L(n,l)T_{2}(l,k)c_{k} \Longleftrightarrow c_{n}=\sum_{k=0}^{n}(-1)^{k}\sum_{l=k}^{n}T_{1}(n,l)L(l,k)a_{k},\\
&a_{n}=\sum_{k=n}^{m}\sum_{l=n}^{k}(-1)^{l}L(k,l)T_{2}(l,n)c_{k} \Longleftrightarrow c_{n}=\sum_{k=n}^{m}(-1)^{n}\sum_{l=n}^{k}T_{1}(k,l)L(l,n)a_{k},
\end{align*}
where $L(n,k)$ are the (unsigned) Lah numbers.                                                                                                  
We would like to remark that in this investigation of central factorial numbers we are led to discover natural definitions for central logarithm and central exponential associated to any delta series. 
The novelty of this paper is that it is the first paper which studies the central factorial numbers of both kinds associated with any sequence of polynomials in a unified and systematic way with the help of umbral calculus. \\
In [26], Koutras generalized the classical Eulerian numbers and polynomials to $p_{n}$-associated Eulerian numbers and polynomials with motivation of providing a unified approach to the study of Eulerian-related numbers and with combinatorial, probabilistic and statistical applications in mind (see [2,9,10,29,33]). In addition, he mentioned in the same paper that a lot of special numbers, like the Stirling numbers of the second, the Lah numbers and Gould-Hopper numbers, appear in the expansions of some polynomials in terms of the falling factorials. This observation is the motivation and impetus for the present research on the central factorial numbers of both kinds associated with sequences of polynomials (see [26,27]). \\
The outline of this paper is as follows. In Section 2, we will briefly go over umbral calculus. In Section 3, we introduce the central factorial numbers of the second associated with sequences of polynomials. We illustrate our results with twenty examples in Section 4. In Section 5, we introduce the central factorial numbers of the first kind associated with sequences of polynomials. We illustrate our results in Section 6 with the same examples from Section 4 and thereby get interesting orthogonality and inverse relations for each case. Finally, we conclude our paper in Section 7.

\section{Review of umbral calculus}
\vspace{0.5cm}
Here we will briefly go over very basic facts about umbral calculus. For more details on this, we recommend the reader to refer to [1,8,31-32]. We remark that recently umbral calculus has been extended to the case of degenerate umbral calculus in order to treat degenerate special polynomials and numbers, which involves degenerate exponentials (see [12]).
Let $\mathbb{C}$ be the field of complex numbers. Then $\mathcal{F}$ denotes the algebra of formal power series in $t$ over $\mathbb{C}$, given by
\begin{displaymath}
 \mathcal{F}=\bigg\{f(t)=\sum_{k=0}^{\infty}a_{k}\frac{t^{k}}{k!}~\bigg|~a_{k}\in\mathbb{C}\bigg\},
\end{displaymath}
and $\mathbb{P}=\mathbb{C}[x]$ indicates the algebra of polynomials in $x$ with coefficients in $\mathbb{C}$. \par
Let $\mathbb{P}^{*}$ be the vector space of all linear functionals on $\mathbb{P}$. If $\langle L|p(x)\rangle$ denotes the action of the linear functional $L$ on the polynomial $p(x)$, then the vector space operations on $\mathbb{P}^{*}$ are defined by
\begin{displaymath}
\langle L+M|p(x)\rangle=\langle L|p(x)\rangle+\langle M|p(x)\rangle,\quad\langle cL|p(x)\rangle=c\langle L|p(x)\rangle,
\end{displaymath}
where $c$ is a complex number. \par
For $f(t)\in\mathcal{F}$ with $\displaystyle f(t)=\sum_{k=0}^{\infty}a_{k}\frac{t^{k}}{k!}\displaystyle$, we define the linear functional on $\mathbb{P}$ by
\begin{equation}\label{1B}
\langle f(t)|x^{k}\rangle=a_{k}. 
\end{equation}
From \eqref{1B}, we note that
\begin{equation*}
 \langle t^{k}|x^{n}\rangle=n!\delta_{n,k},\quad(n,k\ge 0), 
\end{equation*}
where $\delta_{n,k}$ is the Kronecker's symbol. \par
Some remarkable linear functionals are as follows:
\begin{align}
&\langle e^{yt}|p(x) \rangle=p(y), \nonumber \\
&\langle e^{yt}-1|p(x) \rangle=p(y)-p(0), \label{2B} \\
& \bigg\langle \frac{e^{yt}-1}{t}\bigg |p(x) \bigg\rangle = \int_{0}^{y}p(u) du.\nonumber
\end{align}
Let
\begin{equation}\label{3B}
 f_{L}(t)=\sum_{k=0}^{\infty}\langle L|x^{k}\rangle\frac{t^{k}}{k!}.
\end{equation}
Then, by \eqref{1B} and \eqref{3B}, we get
\begin{displaymath}
 \langle f_{L}(t)|x^{n}\rangle=\langle L|x^{n}\rangle.
\end{displaymath}
That is, $f_{L}(t)=L$, as linear functionals on $\mathbb{P}$. In fact, the map $L\longmapsto f_{L}(t)$ is a vector space isomorphism from $\mathbb{P}^{*}$ onto $\mathcal{F}$.\par  Henceforth, $\mathcal{F}$ denotes both the algebra of formal power series  in $t$ and the vector space of all linear functionals on $\mathbb{P}$. $\mathcal{F}$ is called the umbral algebra and the umbral calculus is the study of umbral algebra. 
For each nonnegative integer $k$, the differential operator $t^k$ on $\mathbb{P}$ is defined by
\begin{equation}\label{4B}
t^{k}x^n=\left\{\begin{array}{cc}
(n)_{k}x^{n-k}, & \textrm{if $k\le n$,}\\
0, & \textrm{if $k>n$.}
\end{array}\right. 
\end{equation}
Extending \eqref{4B} linearly, any power series
\begin{displaymath}
 f(t)=\sum_{k=0}^{\infty}\frac{a_{k}}{k!}t^{k}\in\mathcal{F}
\end{displaymath}
gives the differential operator on $\mathbb{P}$ defined by
\begin{equation}\label{5B}
 f(t)x^n=\sum_{k=0}^{n}\binom{n}{k}a_{k}x^{n-k},\quad(n\ge 0). 
\end{equation}
It should be observed that, for any formal power series $f(t)$ and any polynomial $p(x)$, we have
\begin{equation}\label{6B}
\langle f(t) | p(x) \rangle =\langle 1 | f(t)p(x) \rangle =f(t)p(x)|_{x=0}.
\end{equation}
Here we note that an element $f(t)$ of $\mathcal{F}$ is a formal power series, a linear functional and a differential  operator. Some notable differential operators are as follows: 
\begin{align}
&e^{yt}p(x)=p(x+y), \nonumber\\
&(e^{yt}-1)p(x)=p(x+y)-p(x), \label{7B}\\
&\frac{e^{yt}-1}{t}p(x)=\int_{x}^{x+y}p(u) du.\nonumber
\end{align}

The order $o(f(t))$ of the power series $f(t)(\ne 0)$ is the smallest integer for which $a_{k}$ does not vanish. If $o(f(t))=0$, then $f(t)$ is called an invertible series. If $o(f(t))=1$, then $f(t)$ is called a delta series. \par
For $f(t),g(t)\in\mathcal{F}$ with $o(f(t))=1$ and $o(g(t))=0$, there exists a unique sequence $s_{n}(x)$ (deg\,$s_{n}(x)=n$) of polynomials such that
\begin{equation} \label{8B}
\big\langle g(t)f(t)^{k}|s_{n}(x)\big\rangle=n!\delta_{n,k},\quad(n,k\ge 0).
\end{equation}
The sequence $s_{n}(x)$ is said to be the Sheffer sequence for $(g(t),f(t))$, which is denoted by $s_{n}(x)\sim (g(t),f(t))$. We observe from \eqref{8B} that 
\begin{equation}\label{9B}
s_{n}(x)=\frac{1}{g(t)}q_{n}(x),
\end{equation}
where $q_{n}(x)=g(t)s_{n}(x) \sim (1,f(t))$.\par
In particular, if $s_{n}(x) \sim (g(t),t)$, then $q_{n}(x)=x^n$, and hence 
\begin{equation}\label{10B}
s_{n}(x)=\frac{1}{g(t)}x^n.
\end{equation}

It is well known that $s_{n}(x)\sim (g(t),f(t))$ if and only if
\begin{equation}\label{11B}
\frac{1}{g\big(\overline{f}(t)\big)}e^{x\overline{f}(t)}=\sum_{k=0}^{\infty}\frac{s_{k}(x)}{k!}t^{k}, 
\end{equation}
for all $x\in\mathbb{C}$, where $\overline{f}(t)$ is the compositional inverse of $f(t)$ such that $\overline{f}(f(t))=f(\overline{f}(t))=t$. \par

The following equations \eqref{12B}, \eqref{13B}, and \eqref{14B} are equivalent to the fact that  $s_{n}\left(x\right)$ is Sheffer for $\left(g\left(t\right),f\left(t\right)\right)$, for some invertible $g(t)$: 
\begin{align}
f\left(t\right)s_{n}\left(x\right)&=ns_{n-1}\left(x\right),\quad\left(n\ge0\right),\label{12B}\\
s_{n}\left(x+y\right)&=\sum_{j=0}^{n}\binom{n}{j}s_{j}\left(x\right)q_{n-j}\left(y\right),\label{13B}
\end{align}
with $q_{n}\left(x\right)=g\left(t\right)s_{n}\left(x\right),$
\begin{equation}\label{14B}
s_{n}\left(x\right)=\sum_{j=0}^{n}\frac{1}{j!}\big\langle{g\left(\overline{f}\left(t\right)\right)^{-1}\overline{f}\left(t\right)^{j}}\big |{x^{n}\big\rangle}x^{j}.
\end{equation}
If $s_{n}(x)\sim(g(t),f(t))$, then the following recurrence relation holds:
\begin{equation}\label{15B}
s_{n+1}(x)=\bigg(x-\frac{g'(t)}{g(t)}\bigg)\frac{1}{f'(t)}s_{n}(x).
\end{equation}
Let $s_n(x) \sim (g(t),f(t))$. Then, for any polynomial $p(x)$, we have the polynomial exapnsion given by
\begin{align}
p(x)=\sum_{k \ge 0}\frac{1}{k!}\big\langle{g(t)f(t)^{k}|p(x)\big\rangle}s_{k}(x).\label{16B}
\end{align}
and, for any formal power series $h(t)$, we have
\begin{align}
h(t)=\sum_{k=0}^{\infty}\frac{1}{k!}\langle{h(t) | s_{k} \rangle}g(t)f(t)^{k}.\label{17B}
\end{align}
For any formal power series $h(t)$ and any polynomial $p(x)$, we have
\begin{align}
\langle{h(at) | p(x) \rangle}=\langle{h(t) | p(ax) \rangle}.\label{18B}
\end{align}

The falling factorials $(x)_{n}$ are defined by
\begin{align}
(x)_{n}=x(x-1)\cdots(x-n+1), \,\,(n \ge 1),\quad(x)_{0}=1.\label{19B}
\end{align} 
More generally, for any real number $\lambda$, the generalized falling factorials are defined by
\begin{align}
(x)_{n,\lambda}=x(x-\lambda)\cdots(x-(n-1)\lambda), \,\,(n \ge 1),\quad(x)_{0,\lambda}=1\label{20B}
\end{align}
For any nonzero real number $\lambda$, the degenerate exponentials $e_{\lambda}^{x}(t)$ and the degenerate logarithms $\log _{\lambda}t$ are defined by

\begin{align}
& e_{\lambda}^{x}(t)=(1+\lambda t)^{\frac{x}{\lambda}}=\sum_{n=0}^{\infty}(x)_{n,\lambda}\frac{t^{n}}{n!}, \quad e_{\lambda}(t)=e_{\lambda}^{1}(t)=\sum_{n=0}^{\infty}(1)_{n,\lambda}\frac{t^{n}}{n!}, \label{21B}\\
& \log_{\lambda}t=\frac{1}{\lambda}(t^{\lambda}-1). \label{22B}
\end{align}  

We recall that the Stirling numbers of the first kind $S_1(n,k)$ and of the second kind $S_2(n,k)$ are respectively given by 
\begin{align}
&\frac{1}{k!}\big(\log (1+t)\big)^{k}=\sum_{n=k}^{\infty}S_{1}(n,k)\frac{t^{n}}{n!},\quad (x)_{n}=\sum_{k=0}^{n}S_{1}(n,k)x^{k}, \label{23B} \\
&\frac{1}{k!}\big(e^{t}-1\big)^{k}=\sum_{n=k}^{\infty}S_{2}(n,k)\frac{t^{n}}{n!}, \quad x^{n}=\sum_{k=0}^{n}S_2(n,k)(x)_k . \label{24B}
\end{align}

As degenerate versions of the Stirling numbers, the degenerate Stirling numbers of the first kind $S_{1,\lambda}(n,k)$ and of the second kind $S_{2,\lambda}(n,k)$ are respectively given by
\begin{align}
&\frac{1}{k!}\big(\log_{\lambda} (1+t)\big)^{k}=\sum_{n=k}^{\infty}S_{1,\lambda}(n,k)\frac{t^{n}}{n!},\quad (x)_{n}=\sum_{k=0}^{n}S_{1,\lambda}(n,k)(x)_{k,\lambda}, \label{25B} \\
&\frac{1}{k!}\big(e_{\lambda}(t)-1\big)^{k}=\sum_{n=k}^{\infty}S_{2,\lambda}(n,k)\frac{t^{n}}{n!}, \quad (x)_{n,\lambda}=\sum_{k=0}^{n}S_{2,\lambda}(n,k)(x)_k. \label{26B}
\end{align}

The central factorials $x^{[n]}$ are defined by
\begin{align}
x^{[n]}=x(x+\frac{1}{2}n-1)_{n-1}, \,\,(n \ge 1),\quad x^{[0]}=1.\label{27B}
\end{align}
Then, we note that (see [30,31])
\begin{align}
x^{[n]} \sim (1, f(t)=e^{\frac{t}{2}}-e^{-\frac{t}{2}}). \label{28B}
\end{align}
The following facts are well known (see [30,31]). 
\begin{proposition}
Let $f_{1}(t)=e^{\frac{t}{2}}-e^{-\frac{t}{2}}$, and let $f_{2}(t)=e_{\lambda}^{\frac{1}{2}}(t)-e_{\lambda}^{-\frac{1}{2}}(t)$. Then the following hold true.
\begin{align*}
&(a)\,\, \bar{f_{1}}(t)=2\log\bigg(\frac{t+\sqrt{t^{2}+4}}{2}\bigg)=\log \big(1+\frac{t}{2}(t+\sqrt{t^{2}+4}\big)\big). \\
&(b)\,\, \bar{f_{2}}(t)=\log_{\lambda}\bigg(\frac{t+\sqrt{t^{2}+4}}{2}\bigg)^{2}=\frac{1}{\lambda}\bigg(\bigg(\frac{t+\sqrt{t^{2}+4}}{2}\bigg)^{2 \lambda}-1\bigg).\\
&(c)\,\,\sum_{n=0}^{\infty}x^{[n]}\frac{t^{n}}{n!}=\bigg(\frac{t+\sqrt{t^{2}+4}}{2}\bigg)^{2x}=\big(1+\frac{t}{2}(t+\sqrt{t^{2}+4}\big)\big)^{x}.\\
&(d)\,\, e^{y t}=\sum_{k=0}^{\infty}\frac{1}{k!}(e^{\frac{t}{2}}-e^{-\frac{t}{2}})^{k}y^{[k]}.
\end{align*}
\begin{proof}
(a), (b) These can be easily checked. \\
(c) By \eqref{28B}, $\sum_{n=0}^{\infty}x^{[n]}\frac{t^{n}}{n!}=e^{x\bar{f}(t)}$, with $\bar{f}(t)$ as in (a). \\
(d) This follows from \eqref{17B} and \eqref{28B}.
\end{proof}
\end{proposition}

\section{Central factorial  numbers of the second kind associated with sequences of polynomials}

The central factorial numbers of the first kind $T_{1}(n,k)$ and of second kind $T_{2}(n,k)$ are respectively defined by
\begin{align}
&\frac{1}{k!}\bigg(2\log\bigg(\frac{t+\sqrt{t^{2}+4}}{2}\bigg)\bigg)^{k}=\frac{1}{k!}\Big(\log \big(1+\frac{t}{2}(t+\sqrt{t^{2}+4}\big)\big)\Big)^{k} \label{1C}\\
&=\sum_{n=k}^{\infty}T_{1}(n,k)\frac{t^{n}}{n!},\quad x^{[n]}=\sum_{k=0}^{n}T_{1}(n,k)x^{k},\nonumber\\
& \frac{1}{k!}\big(e^{\frac{t}{2}}-e^{-\frac{t}{2}}\big)^{k}=\sum_{n=k}^{\infty}T_{2}(n,k)\frac{t^{n}}{n!},\quad x^{n}=\sum_{k=0}^{n}T_{2}(n,k)x^{[k]}.\label{2C}
\end{align}
We recommend [3-6,30,31] for general references on central factorial numbers. For recent works on central factorial numbers, one refers to [11-13,19-23].
As degenerate versions of the central factorial numbers, the degenerate central factorial numbers of the first kind $T_{1,\lambda}(n,k)$ and of the second kind $T_{2,\lambda}(n,k)$ are respectively defined by
\begin{align}
&\frac{1}{k!}\left\{\frac{1}{\lambda}\bigg(\bigg(\frac{t+\sqrt{t^{2}+4}}{2}\bigg)^{2 \lambda}-1\bigg)\right\}^{k}=\frac{1}{k!}\bigg(\log_{\lambda}\bigg(\frac{t+\sqrt{t^{2}+4}}{2}\bigg)^{2}\bigg)^{k}\label{3C}\\
&\quad=\sum_{n=k}^{\infty}T_{1,\lambda}(n,k)\frac{t^{n}}{n!},
\quad x^{[n]}=\sum_{k=0}^{n}T_{1,\lambda}(n,k)(x)_{k, \lambda},\nonumber\\
&\frac{1}{k!}\big(e_{\lambda}^{\frac{1}{2}}(t)-e_{\lambda}^{-\frac{1}{2}}(t)\big)^{k}=\sum_{n=k}^{\infty}T_{2,\lambda}(n,k)\frac{t^{n}}{n!},\quad (x)_{n,\lambda}=\sum_{k=0}^{n}T_{2,\lambda}(n,k)x^{[k]}.\label{4C}
\end{align}
Let $\bold{P}=\left\{p_{n}(x)\right\}_{n=0}^{\infty}$ be a sequence of polynomials such that deg\,$p_{n}(x)=n, p_{0}(x)=1$.
In \eqref{2C} and \eqref{4C}, we observe that $T_{2}(n,k)$ and $T_{2,\lambda}(n,k)$ arise as the coefficients respectively when we expand $x^{n}$ and $(x)_{n,\lambda}$ in terms of $x^{[k]}$. In view of this observation, it seems natural to define the central factorial numbers of the second associated with $\bold{P}=\left\{p_{n}(x)\right\}_{n=0}^{\infty}$ as the coefficients when we expand $p_{n}(x)$ in terms of $x^{[k]}$:
\begin{align}
p_{n}(x)=\sum_{k=0}^{n}T_2(n,k;\bold{P})x^{[k]}.\label{5C}
\end{align}

\begin{theorem}
Let $\bold{P}=\left\{p_{n}(x)\right\}_{n=0}^{\infty}$ be a sequence of polynomials such that deg\,$p_{n}(x)=n, p_{0}(x)=1$, with $p_{n}(x)=\sum_{l=0}^{n}p_{n,l}x^{l}$. We let
\begin{align*}
p_{n}(x)=\sum_{k=0}^{n}T_2(n,k;\bold{P})x^{[k]}.
\end{align*}
(a)\,\, Then the central factorial numbers of the second kind associated with $\bold{P}$ are given by
\begin{align*}
T_2(n,k;\bold{P})=\frac{1}{k!}\big\langle{(e^{\frac{t}{2}}-e^{-\frac{t}{2}})^{k}|p_{n}(x)\big\rangle}.
\end{align*}
More explicitly, it is given by
\begin{align*}
T_2(n,k;\bold{P})=\sum_{l=k}^{n}T_{2}(l,k)p_{n,l}.
\end{align*}
(b)\,\, $p_{n,l}=\sum_{k=l}^{n}T_{1}(k,l)T_{2}(n,k;\bold{P})$. \\
(c)\,\, $T_2(n,k;\bold{P})$ are alternatively given by
\begin{align*}
T_2(n,k;\bold{P})=\sum_{l=k}^{n}S_{2}(l,k)\frac{1}{l!}\Big(\frac{d}{dx}\Big)^{l}p_{n}(x)\big|_{x=-\frac{k}{2}}.
\end{align*}
(d)\,\,Let $p_{n}(x)$ be Sheffer for the pair $(g(t), f(t))$. Then the generating function of $T_2(n,k;\bold{P})$ is given by
\begin{align*}
\sum_{n=k}^{\infty}T_2(n,k;\bold{P})\frac{t^{n}}{n!}=\frac{1}{g(\bar{f}(t))}\frac{1}{k!}\big(e^{\frac{1}{2}\bar{f}(t)}-e^{-\frac{1}{2}\bar{f}(t)}\big)^{k}=\frac{1}{g(\bar{f}(t))} 
\sum_{n=k}^{\infty}T_{2}(n,k)\frac{\bar{f}(t)^{n}}{n!}.
\end{align*}
(e),\, \textrm{Let $\bar{\bold{P}}=\left\{\bar{p}_{n}(x)\right\}$, where $\bar{p}_{n}(x)$ is the sequence of polynomials defined by} $\bar{p}_{0}(x)=1,\,\,\bar{p}_{n}(x)=xp_{n-1}(x),\,\,(n \ge 1)$. Then we have
\begin{align*}
T_{2}(n+1,k;\bar{\bold{P}})=T_{2}(n,k-1;\bold{P})+\frac{k}{2}T_{2}(n,k;\bold{P}),\,\,(0 \le k \le n+1).
\end{align*}
\begin{proof}
(a)\,\, By noting $x^{[n]} \sim (1, e^{\frac{t}{2}}-e^{-\frac{t}{2}})$, this follows from \eqref{16B}. \\
From \eqref{6B}, \eqref{16B}, and \eqref{2C}, we have
\begin{align}
T_{2}(n,k;\bold{P})&=\frac{1}{k!}\big\langle{(e^{\frac{t}{2}}-e^{-\frac{t}{2}})^{k}|p_{n}(x)\big\rangle}\nonumber\\
&=\sum_{l=k}^{n}T_{2}(l,k)\frac{1}{l!}\big(\frac{d}{dx}\big)^{l}p_{n}(x)|_{x=0}\label{8C}\\
&=\sum_{l=k}^{n}T_{2}(l,k)p_{n,l}.\nonumber
\end{align}
(b) The identity follows from the next observation:
\begin{align*}
\sum_{l=0}^{n}p_{n,l}x^{l}&=\sum_{k=0}^{n}T_{2}(n,k;\bold{P})x^{[k]}\\
&=\sum_{k=0}^{n}T_{2}(n,k;\bold{P})\sum_{l=0}^{k}T_{1}(k,l)x^{l}\\
&=\sum_{l=0}^{n}\sum_{k=l}^{n}T_{1}(k,l)T_{2}(n,k;\bold{P})x^{l}.
\end{align*}
Note here that the explicit expression of $T_{2}(n,k;\bold{P})$ also follows from this by inversion. \\
(c)\,\,From (a) and \eqref{2B}, we have
\begin{align*}
T_2(n,k;\bold{P})&=\frac{1}{k!}\big\langle{(e^{\frac{t}{2}}-e^{-\frac{t}{2}})^{k}|p_{n}(x)\big\rangle}=\frac{1}{k!}\big\langle{(e^{t}-1)^{k}e^{-\frac{k}{2}t}\big| p_{n}(x)\rangle} \\
&=\sum_{l=k}^{n}S_{2}(l,k)\frac{1}{l!}\langle{e^{-\frac{k}{2} t}|t^{l}p_{n}(x)\rangle}\\
&=\sum_{l=k}^{n}S_{2}(l,k)\frac{1}{l!}\Big(\frac{d}{dx}\Big)^{l}p_{n}(x)\big|_{x=-\frac{k}{2}}.
\end{align*}
(d)\,\, From \eqref{11B} and \eqref{5C}, and recalling that $e^{yt}=\sum_{k=0}^{\infty}\frac{1}{k!}(e^{\frac{1}{2}t}-e^{-\frac{1}{2}t})^{k}y^{[k]}$ from Proposition 2.1 (d), we have
\begin{align}
\sum_{k=0}^{\infty}\sum_{n=k}^{\infty}T_2(n,k;\bold{P})\frac{t^{n}}{n!}u^{[k]}&=\sum_{n=0}^{\infty}\sum_{k=0}^{n}T_{2}(n,k;\bold{P})u^{[k]}\frac{t^{n}}{n!}\nonumber\\
&=\sum_{n=0}^{\infty}p_{n}(u)\frac{t^{n}}{n!}=\frac{1}{g(\bar{f}(t))}e^{u\bar{f}(t)}\label{6C}\\
&=\sum_{k=0}^{\infty}\frac{1}{g(\bar{f}(t))}\frac{1}{k!}(e^{\frac{1}{2}\bar{f}(t)}-e^{-\frac{1}{2}\bar{f}(t)})^{k}u^{[k]}.\nonumber
\end{align}
(e)\,\, This follows from the next observation:
\begin{align*}
\sum_{k=0}^{n+1}T_{2}(n+1,k;\bar{\bold{P}})x^{[k]}&=xp_{n}(x)=\sum_{k=0}^{n}(x-\frac{k}{2}+\frac{k}{2})T_{2}(n,k;\bold{P})x^{[k]}\\
&=\sum_{k=0}^{n}T_{2}(n,k;\bold{P})x^{[k+1]}+\sum_{k=0}^{n}\frac{k}{2}T_{2}(n,k;\bold{P})x^{[k]}\\
&=\sum_{k=0}^{n+1}T_{2}(n,k-1;\bold{P})x^{[k]}+\sum_{k=0}^{n+1}\frac{k}{2}T_{2}(n,k;\bold{P})x^{[k]}.
\end{align*} 
\end{proof}
\end{theorem}

\section{Examples on central factorial numbers of the second kind associated with sequences of polynomials}
(a) Let $\bold{P}=\left\{x^{n}\right\}$. Then $x^{n} \sim (1,t)$.
By the definition in \eqref{5C} and Theorem 3.1, we have (see [7,30,31])
\begin{align*}
T_{2}(n,k;\bold{P})=T_{2}(n,k), \quad \sum_{n=k}^{\infty}T_{2}(n,k)\frac{t^n}{n!}=\frac{1}{k!}(e^{\frac{t}{2}}-e^{-\frac{t}{2}})^{k},
\end{align*}
(see \eqref{2C}). \\
(b) Let $\bold{P}=\left\{(x)_{n,\lambda}\right\}$ be the sequence of generalized falling factorials. Then $(x)_{n,\lambda} \sim \big(1, f(t)=\frac{1}{\lambda}(e^{\lambda t}-1)\big)$, with $\bar{f}(t)=\frac{1}{\lambda}\log(1+\lambda t)$ (see [13,14,16,18,23-25]).
By Theorem 3.1, the generating function is given by
\begin{align*}
\sum_{n=k}^{\infty}T_{2}(n,k;\bold{P})\frac{t^{n}}{n!}=\frac{1}{k!}\big(e_{\lambda}^{\frac{1}{2}}(t)-e_{\lambda}^{-\frac{1}{2}}(t)\big)=\sum_{n=k}^{\infty}T_{2,\lambda}(n,k)\frac{t^{n}}{n!},
\end{align*}
(see \eqref{4C}).
Thus $T_{2}(n,k;\bold{P})=T_{2,\lambda}(n,k)$. Alternatively, by Theorem 3.1 it is also given by
\begin{align*}
\sum_{n=k}^{\infty}T_{2}(n,k;\bold{P})\frac{t^{n}}{n!}&=\sum_{l=k}^{\infty}T_{2}(l,k)\frac{1}{l!}\big(\frac{1}{\lambda}\log(1+\lambda t)\big)^{l}\\
&=\sum_{n=k}^{\infty}\sum_{l=k}^{n}\lambda^{n-l}S_{1}(n,l)T_{2}(l,k)\frac{t^{n}}{n!}.
\end{align*}
Thus $T_{2}(n,k;\bold{P})=\sum_{l=k}^{n}\lambda^{n-l}S_{1}(n,l)T_{2}(l,k)$. \\
(c) The rising factorials $\langle x \rangle_{n}$ are defined by (see [7,15,30,31])
\begin{align*}
\langle x \rangle_{n}=x(x+1)\cdots(x+n-1), \,\,(n \ge 1), \quad \langle x \rangle_{0}=1.
\end{align*}
Let $\bold{P}=\left\{\langle{x \rangle}_{n}\right\}$. Then $\langle{x \rangle}_{n} \sim (1, 1-e^{-t}),$ with $\bar{f}(t)=-\log(1-t)$.
From Theorem 3.1, the generating function is given by
\begin{align*}
\sum_{n=k}^{\infty}T_{2}(n,k;\bold{P})\frac{t^{n}}{n!}=\frac{1}{k!}\big((1-t)^{-\frac{1}{2}}-(1-t)^{\frac{1}{2}}\big)^{k}.
\end{align*}
More explicitly, from Theorem 3.1, it is given by
\begin{align*}
\sum_{n=k}^{\infty}T_{2}(n,k;\bold{P})\frac{t^{n}}{n!}&=\sum_{l=k}^{\infty}T_{2}(l,k)\frac{1}{l!}\big(-\log(1-t)\big)^{l}\\
&=\sum_{n=k}^{\infty}\sum_{l=k}^{n}(-1)^{n-l}S_{1}(n,l)T_{2}(l,k)\frac{t^{n}}{n!}.
\end{align*}
Hence we have
\begin{align*}
T_{2}(n,k;\bold{P})=\sum_{l=k}^{n}(-1)^{n-l}S_{1}(n,l)T_{2}(l,k).
\end{align*}
(d) The generalized rising factorials are defined by (see [18])
\begin{align*}
\langle x \rangle_{n,\lambda}=x(x+\lambda)\cdots(x+(n-1)\lambda), \,\,(n \ge 1), \quad \langle x \rangle_{0,\lambda}=1.
\end{align*}
Let $\bold{P}=\left\{\langle{x \rangle}_{n,\lambda}\right\}$. Then $\langle{x \rangle}_{n,\lambda} \sim \big(1, f(t)=\frac{1}{\lambda}(1-e^{-\lambda t})\big)$, with $\bar{f}(t)=-\frac{1}{\lambda}\log(1-\lambda t)$. 
From Theorem 3.1, the generating function is given by
\begin{align*}
\sum_{n=k}^{\infty}T_{2}(n,k;\bold{P})\frac{t^{n}}{n!}=\frac{1}{k!}\Big((1-\lambda t)^{-\frac{1}{2 \lambda}}-(1-\lambda t)^{\frac{1}{2 \lambda}}\Big)^{k}.
\end{align*}
In terms of central factorial numbers of the second, from Theorem 3.1, it  is also given by
\begin{align*}
\sum_{n=k}^{\infty}T_{2}(n,k;\bold{P})\frac{t^{n}}{n!}&=\sum_{l=k}^{\infty}T_{2}(l,k)\frac{1}{l!}\big(-\frac{1}{\lambda}\log(1-\lambda t)\big)^{l}\\
&=\sum_{n=k}^{\infty}\sum_{l=k}^{n}(-\lambda)^{n-l}S_{1}(n,l)T_{2}(l,k)\frac{t^{n}}{n!}.
\end{align*}
Hence we have
\begin{align*}
T_{2}(n,k;\bold{P})=\sum_{l=k}^{n}(-\lambda)^{n-l}S_{1}(n,l)T_{2}(l,k).
\end{align*}

(e) Let $\alpha(t)$ be the delta series given by 
\begin{align}
\alpha(t)=2\Big(\frac{1}{2-t}-\frac{1}{2+t}\Big)=\sum_{n=1}^{\infty}2^{2-2n}t^{2n-1}.\label{4D}
\end{align} 
Then its compositional inverse is given by
\begin{align}
\bar{\alpha}(t)=\frac{2}{t}\big(\sqrt{t^{2}+1}-1\big)=2\sum_{n=1}^{\infty}\big(\frac{1}{2}\big)_{n}\frac{1}{n!}t^{2n-1}.\label{2D}
\end{align}
H. K. Kim (see [19]) defined the central Lah numbers, which we call the central Lah numbers of the second kind, by 
\begin{align}
\frac{1}{k!}(\alpha(t))^{k}=\sum_{n=k}^{\infty}L_{2}^{(c)}(n,k)\frac{t^{n}}{n!}.\label{3D}
\end{align}
We also define the central Lah numbers of the first kind by
\begin{align}
\frac{1}{k!}(\bar{\alpha}(t))^{k}=\sum_{n=k}^{\infty}L_{1}^{(c)}(n,k)\frac{t^{n}}{n!}.\label{4D}
\end{align}
The central Lah-Bell polynomials of the second, which is called the central Lah-Bell polynomials (see [19]), and the central Lah-Bell polynomials of the first kind are respectively defined by
\begin{align*}
LB_{n,2}^{(c)}(x)=\sum_{k=0}^{n}L_{2}^{(c)}(n,k)x^{k}, \quad LB_{n,1}^{(c)}(x)=\sum_{k=0}^{n}L_{1}^{(c)}(n,k)x^{k}.
\end{align*}
Then it is easy to see that their generating functions are given by
\begin{align*}
e^{x\alpha(t)}=\sum_{n=0}^{\infty}LB_{n,2}^{(c)}(x)\frac{t^{n}}{n!},\quad e^{x\bar{\alpha}(t)}=\sum_{n=0}^{\infty}LB_{n,1}^{(c)}(x)\frac{t^{n}}{n!}.
\end{align*}
Let $f(t)=e^{\frac{t}{2}}-e^{-\frac{t}{2}}$. Then its compositional inverse is given by $\bar{f}(t)=2\log\Big(\frac{t+\sqrt{t^{2}+4}}{2}\Big)$.
Define the central factorial Lah numbers of the first kind by
\begin{align}
\frac{1}{k!}\big(\bar{f}(\alpha(t))\big)^{k}&=\frac{1}{k!}\Big(2\log\Big(\frac{1}{2-t}-\frac{1}{2+t}+\sqrt{\big(\frac{1}{2-t}-\frac{1}{2+t}\big)^{2}+1}\Big)\Big)^{k} \label{5D} \\
&=\sum_{n=k}^{\infty}TL_{1}(n,k)\frac{t^{n}}{n!}.\nonumber
\end{align}
Now, we observe from \eqref{1C} and \eqref{3D} that 
\begin{align*}
\sum_{n=k}^{\infty}TL_{1}(n,k)\frac{t^{n}}{n!}&=\sum_{l=k}^{\infty}T_{1}(l,k)\frac{1}{l!}(\alpha(t))^{l}=\sum_{l=k}^{\infty}T_{1}(l,k)\sum_{n=l}^{\infty}L_{2}^{(c)}(n,l)\frac{t^{n}}{n!}\\
&=\sum_{n=k}^{\infty}\sum_{l=k}^{n}T_{1}(l,k)L_{2}^{(c)}(n,l)\frac{t^{n}}{n!}.
\end{align*}
Thus we have shown that 
\begin{align*}
TL_{1}(n,k)=\sum_{l=k}^{n}L_{2}^{(c)}(n,l)T_{1}(l,k).
\end{align*}
Then we define the central factorial Lah Bell polynomials of the first kind $TLB_{n,1}(x)$ as
\begin{align*}
TLB_{n,1}(x)=\sum_{k=0}^{n}TL_{1}(n,k)x^{k}.
\end{align*}
Then it is immediate to see from \eqref{5D} that
\begin{align}
\sum_{n=0}^{\infty}TLB_{n,1}(x)\frac{t^{n}}{n!}=e^{x\bar{f}(\alpha(t))}.\label{6D}
\end{align}
Now, we observe that
\begin{align*}
TLB_{n,1}(x)&=\sum_{k=0}^{n}\sum_{l=k}^{n}L_{2}^{(c)}(n,l)T_{1}(l,k)x^{k}\\
&=\sum_{l=0}^{n}L_{2}^{(c)}(n,l)\sum_{k=0}^{l}T_{1}(l,k)x^{k}
=\sum_{l=0}^{n}L_{2}^{(c)}(n,l)x^{[l]}.
\end{align*}
Hence we have shown that $TLB_{n,1}(x)=\sum_{l=0}^{n}L_{2}^{(c)}(n,l)x^{[l]}$, which says that $T_{2}(n,k;\bold{P_{1}})=L_{2}^{(c)}(n,k)$, with $\bold{P_{1}}=\left\{TLB_{n,1}(x)\right\}_{n=0}^{\infty}$. We note that an explicit expression of this polynomial is obtained in Theorem 8 of [19].  \\
To proceed further, we define the central factorial Lah numbers of the second kind as
\begin{align}
\frac{1}{k!}(\bar{\alpha}(f(t)))^{k}&=\frac{1}{k!}\bigg(\frac{2}{e^{\frac{t}{2}}-e^{-\frac{t}{2}}}\Big(\sqrt{\big(e^{\frac{t}{2}}-e^{-\frac{t}{2}}\big)^{2}+1}-1\Big)\bigg)^{k}\label{7D}\\
&=\frac{1}{k!}\bigg(\frac{2}{e^{\frac{t}{2}}-e^{-\frac{t}{2}}}\Big(\sqrt{e^{t}+e^{-t}-1}-1\Big)\bigg)^{k}\nonumber \\
&=\sum_{n=k}^{\infty}TL_{2}(n,k)\frac{t^{n}}{n!}.\nonumber
\end{align}
Now, we observe from \eqref{2C} and \eqref{4D} that
\begin{align*}
\sum_{n=k}^{\infty}TL_{2}(n,k)\frac{t^{n}}{n!}&=\sum_{l=k}^{\infty}L_{1}^{(c)}(l,k)\frac{1}{l!}(f(t))^{l}
=\sum_{l=k}^{\infty}L_{1}^{(c)}(l,k)\sum_{n=l}^{\infty}T_{2}(n,l)\frac{t^{n}}{n!}\\
&=\sum_{n=k}^{\infty}\sum_{l=k}^{n}T_{2}(n,l)L_{1}^{(c)}(l,k)\frac{t^{n}}{n!}.
\end{align*}
Thus we have shown that
\begin{align*}
TL_{2}(n,k)=\sum_{l=k}^{n}T_{2}(n,l)L_{1}^{(c)}(l,k).
\end{align*}
We define the central factorial Lah Bell polynomials of the second kind $TLB_{n,2}(x)$ as
\begin{align*}
TLB_{n,2}(x)=\sum_{k=0}^{n}TL_{2}(n,k)x^{k}.
\end{align*}
Then it is easy to see from \eqref{7D} that
\begin{align}
\sum_{n=0}^{\infty}TLB_{n,2}(x)\frac{t^{n}}{n!}=e^{x\bar{\alpha}(f(t))}.\label{8D}
\end{align}
Now, we let $\bold{P_{2}}=\left\{TLB_{n,2}(x)\right\}_{n=0}^{\infty}$. Then, from Theorem 3.1 we have $T_{2}(n,k;\bold{P_{2}})=\sum_{l=k}^{n}T_{2}(l,k)TL_{2}(n,l)$.

(f) The central Bell polynomials $\mathrm{Bel}_{n}^{(c)}(x)$ and the central Bell numbers $\mathrm{Bel}_{n}^{(c)}$ are respectively defined by (see [22])
\begin{align*}
e^{x(e^{\frac{t}{2}}-e^{-\frac{t}{2}})}=\sum_{n=0}^{\infty}\mathrm{Bel}_{n}^{(c)}(x)\frac{t^{n}}{n!}, \quad \mathrm{Bel}_{n}^{(c)}=\mathrm{Bel}_{n}^{(c)}(1).
\end{align*}
Then we see that 
\begin{align*}
\mathrm{Bel}_{n}^{(c)}(x) \sim \Big(1, f(t)=2\log\bigg(\frac{t+\sqrt{t^{2}+4}}{2}\bigg)\Big), \quad \mathrm{Bel}_{n}^{(c)}(x)=\sum_{k=0}^{n}T_{2}(n,k) x^{k}.
\end{align*}
Let $\bold{P}=\left\{\mathrm{Bel}_{n}^{(c)}(x)\right\}$. Then, by Theorem 3.1, we see that
\begin{align*}
T_{2}(n,k;\bold{P})=\sum_{l=k}^{n}T_{2}(n,l)T_{2}(l,k).
\end{align*}
Note here that the generating function is given by
\begin{align*}
\sum_{n=k}^{\infty}T_{2}(n,k;\bold{P})\frac{t^{n}}{n!}=\frac{1}{k!}\Big(e^{\frac{1}{2}(e^{\frac{t}{2}}-e^{-\frac{t}{2}})}-e^{-\frac{1}{2}(e^{\frac{t}{2}}-e^{-\frac{t}{2}})}\Big)^{k}.
\end{align*}
(g) The numbers $T_{1,\lambda}(n,k)$ and $T_{2,\lambda}(n,k)$ are respectively called the degenerate central factorial numbers of the first kind and of the second kind, and defined by \eqref{3C} and \eqref{4C} (see [23]).
The degenerate central Bell polynomials $\mathrm{Bel}_{n,\lambda}^{(c)}(x)$ and the degenerate central Bell numbers $\mathrm{Bel}_{n,\lambda}^{(c)}$ are respectively defined by (see [23])
\begin{align*}
e^{x(e_{\lambda}^{\frac{1}{2}}(t)-e_{\lambda}^{-\frac{1}{2}}(t))}=\sum_{n=0}^{\infty}\mathrm{Bel}_{n,\lambda}^{(c)}(x)\frac{t^{n}}{n!}, \quad \mathrm{Bel}_{n,\lambda}^{(c)}=\mathrm{Bel}_{n,\lambda}^{(c)}(1).
\end{align*}
Then we see that 
\begin{align*}
\mathrm{Bel}_{n,\lambda}^{(c)}(x) \sim \Big(1, f(t)=\frac{1}{\lambda}\bigg(\bigg(\frac{t+\sqrt{t^{2}+4}}{2}\bigg)^{2 \lambda}-1\bigg)\Big),\quad
\mathrm{Bel}_{n,\lambda}^{(c)}(x)=\sum_{k=0}^{n}T_{2,\lambda}(n,k) x^{k}.
\end{align*}
Let $\bold{P}=\left\{\mathrm{Bel}_{n}^{(c)}(x)\right\}$. Then, by Theorem 3.1, we see that
\begin{align*}
T_{2}(n,k;\bold{P})=\sum_{l=k}^{n}T_{2}(l,k)T_{2,\lambda}(n,l).
\end{align*}
We note here that the generating function is given by
\begin{align*}
\sum_{n=k}^{\infty}T_{2}(n,k;\bold{P})\frac{t^{n}}{n!}=\frac{1}{k!}\Big(e^{\frac{1}{2}\big(e_{\lambda}^{\frac{1}{2}}(t)-e_{\lambda}^{-\frac{1}{2}}(t)\big)}-e^{-\frac{1}{2}\big(e_{\lambda}^{\frac{1}{2}}(t)-e_{\lambda}^{-\frac{1}{2}}(t)\big)} \Big)^{k}.
\end{align*}
(h) Let $x^{[n,\lambda]}$ be defined by
\begin{align*}
x^{[n,\lambda]}=x\big(x+(\frac{1}{2}n-1)\lambda\big)_{n-1,\lambda}, \,\,(n \ge 1),\quad x^{[0,\lambda]}=1.
\end{align*}
where $x^{[n,\lambda]} \sim \big(1, f(t)=\frac{1}{\lambda}(e^{\frac{\lambda t}{2}}-e^{-\frac{\lambda t}{2}})\big)$. Here we note that
\begin{align*}
\bar{f}(t)=\frac{2}{\lambda}\log\bigg(\frac{\lambda t+\sqrt{\lambda^{2}t^{2}+4}}{2}\bigg)=\frac{1}{\lambda}\log \big(1+\frac{\lambda t}{2}(\lambda t+\sqrt{\lambda^{2}t^{2}+4})\big).
\end{align*}
As analogies to the numbers $T_{1,\lambda}(n,k)$ and $T_{2,\lambda}(n,k)$, we may introduce the numbers $R_{1,\lambda}(n,k)$ and $R_{2,\lambda}(n,k)$ which are defined by
\begin{align*}
&x^{[n,\lambda]}=\sum_{k=0}^{n}R_{1,\lambda}(n,k)x^{k},\quad \sum_{n=k}^{\infty}R_{1,\lambda}(n,k)\frac{t^{n}}{n!}=\frac{1}{k!}\bigg(\frac{2}{\lambda}\log\bigg(\frac{\lambda t+\sqrt{\lambda^{2}t^{2}+4}}{2}\bigg)\bigg)^{k},\\
&x^{n}=\sum_{k=0}^{n}R_{2,\lambda}(n,k)x^{[k,\lambda]}, \quad \sum_{n=k}^{\infty}R_{2,\lambda}(n,k)\frac{t^{n}}{n!}=\frac{1}{k!}\Big(\frac{1}{\lambda}\big(e^{\frac{\lambda t}{2}}-e^{-\frac{\lambda t}{2}}\big)\Big)^{k}.
\end{align*}
Let $\bold{P}=\left\{x^{[n,\lambda]}\right\}$. From Theorem 3.1, we get
\begin{align*}
T_{2}(n,k;\bold{P})=\sum_{l=k}^{n}T_{2}(l,k)R_{1,\lambda}(n,l).
\end{align*}
From Theorem 3.1, the generating function is given by
\begin{align*}
\sum_{n=k}^{\infty}T_{2}(n,k;\bold{P})\frac{t^{n}}{n!}=\frac{1}{k!}\bigg(\bigg(\frac{\lambda t+\sqrt{\lambda^{2}t^{2}+4}}{2}\bigg)^{\frac{1}{\lambda}}-\bigg(\frac{\lambda t+\sqrt{\lambda^{2}t^{2}+4}}{2}\bigg)^{-\frac{1}{\lambda}}\bigg)^{k}.
\end{align*}
(i) Let $\bold{P}=\left\{B_{n}^{L}(x)\right\}$ be the sequence of Lah-Bell polynomials. Then it is given by $e^{x(\frac{t}{1-t})}=\sum_{n=0}^{\infty}B_{n}^{L}(x)\frac{t^n}{n!}$, so that $B_{n}^{L}(x) \sim (1,\frac{t}{1+t})$, and $B_{n}^{L}(x)=\sum_{k=0}^{n}L(n,k)x^{k}$. Here $L(n,k)=\binom{n-1}{k-1}\frac{n!}{k!}$ are the (unsigned) Lah numbers (see [7,15,30,31]) whose generating function is given by $\sum_{n=k}^{\infty}L(n,k)\frac{t^{n}}{n!}=\frac{1}{k!}\big(\frac{t}{1-t}\big)^{k}$.  By Theorem 3.1, we have
\begin{align*}
T_{2}(n,k;\bold{P})=\sum_{l=k}^{n}T_{2}(l,k)L(n,l).
\end{align*}
In addition, by Theorem 3.1, the generating function is given by
\begin{align*}
\sum_{n=k}^{\infty}T_{2}(n,k;\bold{P})\frac{t^{n}}{n!}=\frac{1}{k!}(e^{\frac{t}{2(1-t)}}-e^{-\frac{t}{2(1-t)}})^{k}.
\end{align*}
(j) Let $\bold{P}=\left\{B_{n,\lambda}^{L}(x)\right\}$ be the sequence of degenerate Lah-Bell polynomials. Then it is given by $e_{\lambda}^{x}(\frac{t}{1-t})=\sum_{n=0}^{\infty}B_{n,\lambda}^{L}(x)\frac{t^n}{n!}$. Thus we see that $B_{n,\lambda}^{L}(x) \sim (1,f(t)=\frac{e^{\lambda t}-1}{\lambda+e^{\lambda t}-1})$, with $\bar{f}(t)=\frac{1}{\lambda}\log(1+\frac{\lambda t}{1-t})$ (see [18]).  By Theorem 3.1,  we have
\begin{align*}
\sum_{n=k}^{\infty}T_{2}(n,k;\bold{P})\frac{t^{n}}{n!}&=\sum_{l=k}^{\infty}T_{2}(l,k)\frac{1}{l!}\Big(\frac{1}{\lambda}\log\big(1+\frac{\lambda t}{1-t}\big)\Big)^{l}\\
&=\sum_{l=k}^{\infty}T_{2}(l,k)\frac{1}{\lambda^{l}}\sum_{m=l}^{\infty}S_{1}(m,l)\frac{1}{m!}\big(\frac{\lambda t}{1-t}\big)^{m}\\
&=\sum_{l=k}^{\infty}\sum_{m=l}^{\infty}\lambda^{m-l}S_{1}(m,l)T_{2}(l,k)\sum_{n=m}^{\infty}L(n,m)\frac{t^{n}}{n!}\\
&=\sum_{n=k}^{\infty}\sum_{l=k}^{n}\sum_{m=l}^{n}\lambda^{m-l}S_{1}(m,l)T_{2}(l,k)L(n,m)\frac{t^{n}}{n!}.
\end{align*}
Thus we obtain $T_{2}(n,k;\bold{P})=\sum_{l=k}^{n}\sum_{m=l}^{n}\lambda^{m-l}L(n,m)S_{1}(m,l)T_{2}(l,k).$ \\
Alternatively, the generating function is also given by
\begin{align*}
\sum_{n=k}^{\infty}T_{2}(n,k;\bold{P})\frac{t^{n}}{n!}=\frac{1}{k!}\Big(\Big(1+\frac{\lambda t}{1-t}\Big)^{\frac{1}{2 \lambda}}-\Big(1+\frac{\lambda t}{1-t}\Big)^{-\frac{1}{2 \lambda}}\Big)^{k}.
\end{align*}
(k) Let $\bold{P}=\left\{\mathrm{Bel}_{n}(x)\right\}$ be the sequence of Bell polynomials given by $e^{x(e^{t}-1)}=\sum_{n=0}^{\infty}\mathrm{Bel}_{n}(x)\frac{t^{n}}{n!}$. Then $\mathrm{Bel}_{n}(x) \sim (1, \log (1+t))$, and $\mathrm{Bel}_{n}(x)=\sum_{k=0}^{\infty}S_{2}(n,k)x^{k}$ (see [7,30,31]).  By Theorem 3.1, the generating function is given by
\begin{align*}
\sum_{n=k}^{\infty}T_{2}(n,k;\bold{P})\frac{t^{n}}{n!}&=\frac{1}{k!}\big(e^{\frac{1}{2}(e^t-1)}-e^{-\frac{1}{2}(e^t-1)}\big)^{k}.
\end{align*}
More explicitly, from Theorem 3.1, we get
\begin{align*}
T_{2}(n,k;\bold{P})=\sum_{l=k}^{n}S_{2}(n,l)T_{2}(l,k).
\end{align*}
(l) Let $\bold{P_{1}}=\left\{\mathrm{Bel}_{n,\lambda}(x)\right\}$ be the sequence of partially degenerate Bell polynomials given by  $e^{x(e_{\lambda}(t)-1)}=\sum_{n=0}^{\infty}\mathrm{Bel}_{n,\lambda}(x)\frac{t^{n}}{n!}$. Then $\mathrm{Bel}_{n,\lambda}(x) \sim (1, \log_{\lambda}(1+t))$, and $\mathrm{Bel}_{n,\lambda}(x)=\sum_{k=0}^{n}S_{2,\lambda}(n,k)x^{k}$ (see [24]).  From Theorem 3.1, it is immediate to see that
\begin{align*}
T_{2}(n,k;\bold{P_{1}})=\sum_{l=k}^{n}S_{2,\lambda}(n,l)T_{2}(l,k).
\end{align*}
The generating function is given by
\begin{align*}
\sum_{n=k}^{\infty}T_{2}(n,k;\bold{P_{1}})\frac{t^{n}}{n!}=\frac{1}{k!}\big(e^{\frac{1}{2}(e_{\lambda}(t)-1)}-e^{-\frac{1}{2}(e_{\lambda}(t)-1)}\big)^{k}.
\end{align*}
There are several other degenerate versions of Bell polynomials, one of which is the so called fully degenerate Bell polynomials $\phi_{n,\lambda}(x)$ given by $e_{\lambda}^{x}(e_{\lambda}(t)-1)=\sum_{n=0}^{\infty}\phi_{n,\lambda}(x)\frac{t^{n}}{n!}.$ Then $\phi_{n,\lambda}(x) \sim \big(1, f(t)=\log_{\lambda}(1+\frac{1}{\lambda}(e^{\lambda t}-1))\big)$, with $\bar{f}(t)= \frac{1}{\lambda}\log(1+\lambda(e_{\lambda}(t)-1)\big)$
(see [25]). Let $\bold{P_{2}}=\left\{\phi_{n,\lambda}(x)\right\}$. From Theorem 3.1, we have
\begin{align*}
\sum_{n=k}^{\infty}T_{2}(n,k;\bold{P_{2}})\frac{t^{n}}{n!}&=\sum_{l=k}^{\infty}T_{2}(l,k)\frac{1}{l!}\Big(\frac{1}{\lambda}\log\big(1+\lambda(e_{\lambda}(t)-1)\big)\Big)^{l}\\
&=\sum_{l=k}^{\infty}T_{2}(l,k)\frac{1}{\lambda^{l}}\sum_{m=l}^{\infty}S_{1}(m,l)\frac{1}{m!}\big(\lambda(e_{\lambda}(t)-1)\big)^{m}\\
&=\sum_{l=k}^{\infty}\sum_{m=l}^{\infty}T_{2}(l,k)S_{1}(m,l)\lambda^{m-l}\sum_{n=m}^{\infty}S_{2,\lambda}(n,m)\frac{t^{n}}{n!}\\
&=\sum_{n=k}^{\infty}\sum_{l=k}^{n}\sum_{m=l}^{n}\lambda^{m-l}S_{2,\lambda}(n,m)S_{1}(m,l)T_{2}(l,k)\frac{t^{n}}{n!},
\end{align*}
which verifies that we have $T_{2}(n,k;\bold{P_{2}})=\sum_{l=k}^{n}\sum_{m=l}^{n}\lambda^{m-l}S_{2,\lambda}(n,m)S_{1}(m,l)T_{2}(l,k).$\\
Alternatively, the generating function is also given by
\begin{align*}
\sum_{n=k}^{\infty}T_{2}(n,k;\bold{P_{2}})\frac{t^{n}}{n!}=\frac{1}{k!}\Big(\big(1+\lambda(e_{\lambda}(t)-1)\big)^{\frac{1}{2 \lambda}}-\big(1+\lambda(e_{\lambda}(t)-1)\big)^{-\frac{1}{2 \lambda}}\Big)^{k}.
\end{align*}
(m) Let $\bold{P}=\left\{M_{n}(x)\right\}$ be the sequence of Mittag-Leffler polynomials. That is, $M_{n}(x) \sim (1,f(t)=\frac{e^{t}-1}{e^{t}+1})$, with $\bar{f}(t)=\log \big(\frac{1+t}{1-t}\big)$ (see [31]).
By Theorem 3.1, we have
\begin{align*}
\sum_{n=k}^{\infty}T_{2}(n,k;\bold{P})\frac{t^{n}}{n!}=\frac{1}{k!}\Big(\Big(\frac{1+t}{1-t}\Big)^{\frac{1}{2}}-\Big(\frac{1+t}{1-t}\Big)^{-\frac{1}{2}}\Big)^{k}.
\end{align*}
Alternatively, we also have
\begin{align*}
\sum_{n=k}^{\infty}T_{2}(n,k;\bold{P})\frac{t^{n}}{n!}&=\sum_{l=k}^{\infty}T_{2}(l,k)\frac{1}{l!}\Big(\log\Big(1+\frac{2t}{1-t}\Big)\Big)^{l}\\
&=\sum_{l=k}^{\infty}T_{2}(l,k)\sum_{m=l}^{\infty}S_{1}(m,l)2^{m}\frac{1}{m!}\Big(\frac{t}{1-t}\Big)^{m}\\
&=\sum_{l=k}^{\infty}T_{2}(l,k)\sum_{m=l}^{\infty}S_{1}(m,l)2^{m}\sum_{n=m}^{\infty}L(n,m)\frac{t^{n}}{n!}\\
&=\sum_{n=k}^{\infty}\sum_{l=k}^{n}\sum_{m=l}^{n}2^{m}L(n,m)S_{1}(m,l)T_{2}(l,k)\frac{t^{n}}{n!}.
\end{align*}
This shows that $T_{2}(n,k;\bold{P})=\sum_{l=k}^{n}\sum_{m=l}^{n}2^{m}L(n,m)S_{1}(m,l)T_{2}(l,k)$.

(n) Let $\bold{P}=\left\{L_{n}(x)\right\}$ be the sequence of Laguerre polynomials of order -1 (see [31]).
That is, $L_{n}(x) \sim (1, f(t)=\frac{t}{t-1})$, with $\bar{f}(t)=\frac{t}{t-1}$. By Theorem 3.1, we have
\begin{align*}
\sum_{n=k}^{\infty}T_{2}(n,k;\bold{P})\frac{t^{n}}{n!}&=\sum_{l=k}^{\infty}T_{2}(l,k)(-1)^{l}\frac{1}{l!}\Big(\frac{t}{1-t}\Big)^{l}\\
&=\sum_{l=k}^{\infty}T_{2}(l,k)(-1)^{l}\sum_{n=l}^{\infty}L(n,l)\frac{t^{n}}{n!}\\
&=\sum_{n=k}^{\infty}\sum_{l=k}^{n}(-1)^{l}L(n,l)T_{2}(l,k)\frac{t^{n}}{n!}.
\end{align*}
Thus we have $T_{2}(n,k;\bold{P})=\sum_{l=k}^{n}(-1)^{l}L(n,l)T_{2}(l,k)$. Moreover, the generating function is also given by
\begin{align*}
\sum_{n=k}^{\infty}T_{2}(n,k;\bold{P})\frac{t^{n}}{n!}=\frac{1}{k!}\big(e^{\frac{t}{2(t-1)}}-e^{-\frac{t}{2(t-1)}}\big)^{k}.
\end{align*}
(o) Let $\bold{P}=\left\{B_{n}(x)\right\}$ be the sequence of Bernoulli polynomials. Then $B_{n}(x) \sim \big(\frac{e^{t}-1}{t},t\big)$, and $B_{n}(x)=\sum_{k=0}^{n}\binom{n}{k}B_{n-k}x^{k}$ (see [31]).  By Theorem 3.1, we have
\begin{align*}
T_{2}(n,k;\bold{P})&=\sum_{l=k}^{n}T_{2}(l,k)\binom{n}{l}B_{n-l}=\sum_{l=0}^{n-k}\binom{n}{l}T_{2}(n-l,k)B_{l},\\
&\sum_{n=k}^{\infty}T_{2}(n,k;\bold{P})\frac{t^{n}}{n!}=\frac{t}{e^{t}-1}\frac{1}{k!}(e^{\frac{t}{2}}-e^{-\frac{t}{2}})^{k}.
\end{align*}
(p) Let $\bold{P}=\left\{E_{n}(x)\right\}$ be the sequence of Euler polynomials. Then $E_{n}(x) \sim \big(\frac{e^{t}+1}{2},t\big)$ (see [31]).
Analogously to (c), by Theorem 3.1 we see that
\begin{align*}
&T_{2}(n,k;\bold{P})=\sum_{l=k}^{n}T_{2}(l,k)\binom{n}{l}E_{n-l}=\sum_{l=0}^{n-k}\binom{n}{l}T_{2}(n-l,k)E_{l},\\
&\sum_{n=k}^{\infty}T_{2}(n,k;\bold{P})\frac{t^{n}}{n!}=\frac{2}{e^{t}+1}\frac{1}{k!}(e^{\frac{t}{2}}-e^{-\frac{t}{2}})^{k}.
\end{align*}
(q) Let $\bold{P}=\left\{(rx+s)_{n}\right\}$, with $r \ne 0$. Then
\begin{align*}
(rx+s)_{n}=\sum_{k=0}^{n}G(n,k;r,s)(x)_{k}, 
\end{align*}
where $G(n,k;r,s)$ are called the Gould-Hopper numbers (see [26]).  
We note that $(rx+s)_{n} \sim \big(e^{-\frac{s}{r} t}, f(t)=e^{\frac{t}{r}}-1\big)$, where $\bar{f}(t)=r\log(1+t)$. By Theorem 3.1, the generating function is given by
\begin{align*}
\sum_{n=k}^{\infty}T_{2}(n,k;\bold{P})\frac{t^{n}}{n!}=(1+t)^{s}\frac{1}{k!}\big((1+t)^{\frac{r}{2}}-(1+t)^{-\frac{r}{2}}\big)^{k}.
\end{align*}
Alternatively, from Theorem 3.1 the generating function is also given by
\begin{align*}
\sum_{n=k}^{\infty}T_{2}(n,k;\bold{P}) \frac{t^{n}}{n!}&=(1+t)^{s}\sum_{l=k}^{\infty}T_{2}(l,k)r^{l}\frac{1}{l!}(\log(1+t))^{l}\\
&=(1+t)^{s}\sum_{l=k}^{\infty}T_{2}(l,k)r^{l}\sum_{m=l}^{\infty}S_{1}(m,l)\frac{t^{m}}{m!}\\
&=\sum_{j=0}^{\infty}(s)_{j}\frac{t^{j}}{j!}\sum_{m=k}^{\infty}\sum_{l=k}^{m}r^{l}S_{1}(m,l)T_{2}(l,k)\frac{t^{m}}{m!}\\
&=\sum_{n=k}^{\infty}\sum_{l=k}^{n}\sum_{m=l}^{n}\binom{n}{m}r^{l}(s)_{n-m}S_{1}(m,l)T_{2}(l,k)\frac{t^{n}}{n!}.
\end{align*}
Thus we have
\begin{align*}
T_{2}(n,k;\bold{P})=\sum_{l=k}^{n}\sum_{m=l}^{n}\binom{n}{m}r^{l}(s)_{n-m}S_{1}(m,l)T_{2}(l,k).
\end{align*}
(r) Let $\bold{P}=\left\{b_{n}(x)\right\}$ be the sequence of Bernoulli polynomials of the second kind. Then $b_{n}(x) \sim \big(\frac{t}{e^{t}-1}, f(t)=e^{t}-1\big)$, with $\bar{f}(t)=\log(1+t)$, so that $\sum_{n=0}^{\infty}b_{n}(x)\frac{t^{n}}{n!}=\frac{t}{\log (1+t)}(1+t)^{x}$ (see [31]).  From Theorem 3.1, the generating function is
\begin{align*}
\sum_{n=k}^{\infty}T_{2}(n,k;\bold{P})\frac{t^{n}}{n!}=\frac{t}{\log(1+t)}\frac{1}{k!}\big((1+t)^{\frac{1}{2}}-(1+t)^{-\frac{1}{2}}\big)^{k}.
\end{align*}
Alternatively, we also have
\begin{align*}
\sum_{n=k}^{\infty}T_{2}(n,k;\bold{P})\frac{t^{n}}{n!}&=\frac{t}{\log (1+t)}\sum_{l=k}^{\infty}T_{2}(l,k)\frac{1}{l!}(\log(1+t))^{l}\\
&=\frac{t}{\log (1+t)}\sum_{l=k}^{\infty}T_{2}(l,k)\sum_{m=l}^{\infty}S_{1}(m,l)\frac{t^{m}}{m!}\\
&=\sum_{j=0}^{\infty}b_{j}\frac{t^{j}}{j!}\sum_{m=k}^{\infty}\sum_{l=k}^{m}S_{1}(m,l)T_{2}(l,k)\frac{t^{m}}{m!}\\
&=\sum_{n=k}^{\infty}\sum_{l=k}^{n}\sum_{m=l}^{n}\binom{n}{m}b_{n-m}S_{1}(m,l)T_{2}(l,k)\frac{t^{n}}{n!}.
\end{align*}
Thus we get
\begin{align*}
T_{2}(n,k;\bold{P})=\sum_{l=k}^{n}\sum_{m=l}^{n}\binom{n}{m}b_{n-m}S_{1}(m,l)T_{2}(l,k),
\end{align*}
where $b_{n}=b_{n}(0)$ are the Bernoulli numbers of the second. \\
(s) Let $\bold{P}=\left\{C_{n}(x;a)\right\}$ be the sequence of Poisson-Charlier polynomials. Here $C_{n}(x;a) \sim \big(e^{a(e^{t}-1)}, f(t)=a(e^{t}-1)\big)$, with $a \ne 0$. Then $\bar{f}(t)=\log(1+\frac{t}{a})$, and $\sum_{n=0}^{\infty}C_{n}(x;a)\frac{t^{n}}{n!}=e^{-t}\big(1+\frac{t}{a}\big)^{x}$ (see [31]).  From Theorem 3.1, the generating function is given by
\begin{align*}
\sum_{n=k}^{\infty}T_{2}(n,k;\bold{P})\frac{t^{n}}{n!}=e^{-t}\frac{1}{k!}\Big(\Big(1+\frac{t}{a}\Big)^{\frac{1}{2}}-\Big(1+\frac{t}{a}\Big)^{-\frac{1}{2}}\Big)^{k}.
\end{align*}
Another expression for generating function follows from Theorem 3.1, which is given by
\begin{align*}
\sum_{n=k}^{\infty}T_{2}(n,k;\bold{P})\frac{t^{n}}{n!}&=e^{-t}\sum_{l=k}^{\infty}T_{2}(l,k)\frac{1}{l!}\big(\log(1+\frac{t}{a})\big)^{l}\\
&=e^{-t}\sum_{l=k}^{\infty}T_{2}(l,k)\sum_{m=l}^{\infty}S_{1}(m,l)a^{-m}\frac{t^{m}}{m!}\\
&=\sum_{j=0}^{\infty}(-1)^{j}\frac{t^{j}}{j!}\sum_{m=k}^{\infty}\sum_{l=k}^{m}a^{-m}S_{1}(m,l)T_{2}(l,k)\frac{t^{m}}{m!}\\
&=\sum_{n=k}^{\infty}\sum_{l=k}^{n}\sum_{m=l}^{n}\binom{n}{m}(-1)^{n-m}a^{-m}S_{1}(m,l)T_{2}(l,k)\frac{t^{n}}{n!}.
\end{align*}
Thus we see that
\begin{align*}
T_{2}(n,k;\bold{P})=\sum_{l=k}^{n}\sum_{m=l}^{n}\binom{n}{m}(-1)^{n-m}a^{-m}S_{1}(m,l)T_{2}(l,k).
\end{align*}
(t) Let $\bold{P}=\left\{p_{n}(x)\right\}$, with $p_{n}(x)=\sum_{k=0}^{n}B_{k}(x)B_{n-k}(x)$. This is not a Sheffer sequence. For this, we recall from [17] that
\begin{align*}
p_{n}(x)=\frac{2}{n+2}\sum_{m=0}^{n-2}\binom{n+2}{m}B_{n-m}B_m(x)+(n+1)B_n(x).
\end{align*}
By \eqref{12B} and \eqref{2C}, we note that
\begin{align*}
\frac{1}{k!}\langle{(e^{\frac{t}{2}}-e^{-\frac{t}{2}})^{k}|B_{n}(x) \rangle}&=\sum_{l=0}^{n-k}\binom{n}{l}T_{2}(n-l,k)B_{l}.
\end{align*}
By Theorem 3.1, we see that
\begin{align*}
&T_{2}(n,k;\bold{P})=\frac{1}{k!}\langle{(e^{\frac{t}{2}}-e^{-\frac{t}{2}})^{k}|p_{n}(x)\rangle}\\
&=\frac{2}{n+2}\sum_{m=0}^{n-2}\binom{n+2}{m}B_{n-m}\frac{1}{k!}\langle{(e^{\frac{t}{2}}-e^{-\frac{t}{2}})^{k}|B_{m}(x) \rangle}+(n+1)\frac{1}{k!}\langle{(e^{\frac{t}{2}}-e^{-\frac{t}{2}})^{k}|B_{n}(x) \rangle}\\
&=\frac{2}{n+2}\sum_{m=0}^{n-2}\binom{n+2}{m}B_{n-m}\sum_{l=0}^{m-k}\binom{m}{l}T_{2}(m-l,k)B_{l}
+(n+1)\sum_{l=0}^{n-k}\binom{n}{l}T_{2}(n-l,k)B_{l}.
\end{align*}
This says that the following identity holds:
\begin{align*}
\sum_{k=0}^{n}B_{k}(x)B_{n-k}(x)&=\sum_{k=0}^{n}\left\{\frac{2}{n+2}\sum_{m=0}^{n-2}\binom{n+2}{m}B_{n-m}\sum_{l=0}^{m-k}\binom{m}{l}T_{2}(m-l,k)B_{l}\right.\\
&\left.\quad\quad\quad+(n+1)\sum_{l=0}^{n-k}\binom{n}{l}T_{2}(n-l,k)B_{l}\right\}x^{[k]}.
\end{align*}

\section{Central factorial numbers of the first kind associated with sequences of polynomials}

As in the `classical' case, we would like to introduce the central factorial numbers of the first kind associated with sequences of polynomials. Let $\bold{P}=\left\{p_{n}(x)\right\}_{n=0}^{\infty}$ be a sequence of polynomials such that deg\,$p_{n}(x)=n, p_{0}(x)=1$.
In view of our definition of the central factorial numbers of the second kind associated with $\bold{P}=\left\{p_{n}(x)\right\}_{n=0}^{\infty}$, it is natural to define the central factorial numbers of the first kind associated with $\bold{P}=\left\{p_{n}(x)\right\}_{n=0}^{\infty}$ as the coefficients when we expand $x^{[n]}$ in terms of $p_{k}(x)$:
\begin{align}
x^{[n]}=\sum_{k=0}^{n}T_1(n,k;\bold{P})p_{k}(x).\label{1E}
\end{align}

\begin{theorem}
Let $\bold{P}=\left\{p_{n}(x)\right\}_{n=0}^{\infty}$ be a sequence of polynomials such that deg\,$p_{n}(x)=n, p_{0}(x)=1$. We let
\begin{align*}
x^{[n]}=\sum_{k=0}^{n}T_1(n,k;\bold{P})p_{k}(x).
\end{align*}
(a)\,\, Let $\bar{\bold{P}}=\left\{\bar{p}_{n}(x)\right\}$, where $\bar{p}_{n}(x)$ is the sequence of polynomials defined by $\bar{p}_{0}(x)=1,\,\,\bar{p}_{n}(x)=xp_{n-1}(x),\,\,(n \ge 1)$. Then we have
\begin{align*}
T_{1}(n+1,k;\bar{\bold{P}})=T_{1}(n,k-1;\bold{P})-\frac{n}{2}T_{1}(n,k;\bar{\bold{P}}).
\end{align*}
(b)\,\, Assume that $p_{n}(x) \sim (g(t), f(t))$. Then the central factorial numbers of the first kind associated with $\bold{P}$ are given by
\begin{align*}
T_1(n,k;\bold{P})=\frac{1}{k!}\langle{g(t)(f(t))^{k}|x^{[n]}\rangle}=\frac{1}{k!}\sum_{l=k}^{n}T_{1}(n,l)\langle{g(t)f(t)^{k}|x^{l}\rangle}.
\end{align*}
(c)\,\,$\sum_{k=0}^{n}T_{1}(n,k;\bold{P})p_{k}(1)=(\frac{n}{2})_{n-1},\,\,(n \ge 1).$
\begin{proof}
(a)\,\,This follows from the next observation:
\begin{align*}
\sum_{k=0}^{n+1}T_{1}(n+1,k;\bar{\bold{P}})\bar{p}_{k}(x)&=x^{[n]}\big(x-\frac{n}{2}\big)\\
&=x\sum_{k=0}^{n}T_{1}(n,k;\bold{P})p_{k}(x)-\frac{n}{2}\sum_{k=0}^{n}T_{1}(n,k;\bar{\bold{P}})\bar{p}_{k}(x)\\
&=\sum_{k=0}^{n}T_{1}(n,k;\bold{P})\bar{p}_{k+1}(x)-\frac{n}{2}\sum_{k=0}^{n}T_{1}(n,k;\bar{\bold{P}})\bar{p}_{k}(x)\\
&=\sum_{k=0}^{n+1}T_{1}(n,k-1;\bold{P})\bar{p}_{k}(x)-\frac{n}{2}\sum_{k=0}^{n+1}T_{1}(n,k;\bar{\bold{P}})\bar{p}_{k}(x).
\end{align*}
(b)\,\, This follows from \eqref{16B} and \eqref{1C}. \\
(c)\,\, This is clear by definition in \eqref{1E} and \eqref{27B}.
\end{proof}
\end{theorem}
For any delta series $f(t)$, we define the {\it{central logarithm associated to $f(t)$}} by
\begin{align*}
LC_{f(t)} t = f \big(2 \log \Big(\frac{t+\sqrt{t^{2}+4}}{2}\Big)\big),
\end{align*}
In case $p_{n}(x) \sim (1, f(t))$, $LC_{f(t)}t$ may be also denoted by $LC_{\bold{P}}t$ and called the {\it{central logarithm associated to $\bold{P}=\left\{p_{n}(x)\right\}$}}. Recalling that $x^{[n]} \sim (1, l(t)=e^{\frac{t}{2}}-e^{-\frac{t}{2}})$, we remark that $LC_{f(t)}t=f(\bar{l}(t))$.
\begin{theorem}
Let $p_{n}(x) \sim (1, f(t)),$  and let $\bold{P}=\left\{p_{n}(x)\right\}$. Then the generating function of the central factorial numbers of the first kind associated with $\bold{P}$ is given by
\begin{align*}
\sum_{n=k}^{\infty}T_1(n,k;\bold{P})\frac{t^{n}}{n!}=\frac{1}{k!}(LC_{f(t)}t)^{k},
\end{align*}
where $LC_{f(t)}t=LC_{\bold{P}}t=f \big(2 \log \Big(\frac{t+\sqrt{t^{2}+4}}{2}\Big)\big)$ is the central logarithm associated to $f(t)$.
\begin{proof}
Assume now that $p_{n}(x) \sim (1, f(t))$. Then we observe that
\begin{align*}
\sum_{k=0}^{\infty}p_{k}(x)\sum_{n=k}^{\infty}T_1(n,k;\bold{P})\frac{t^{n}}{n!}&=\sum_{n=0}^{\infty}\sum_{k=0}^{n}T_{1}(n,k;\bold{P})p_{k}(x)\frac{t^{n}}{n!}\\
&=\sum_{n=0}^{\infty}x^{[n]}\frac{t^{n}}{n!}=\Big(\frac{t+\sqrt{t^{2}+4}}{2}\Big)^{2 x}.
\end{align*}
On the other hand, we also have
\begin{align*}
\sum_{k=0}^{\infty}p_{k}(x)\frac{1}{k!}\big(LC_{f(t)}t\big)^{k}&=e^{x\bar{f}\big(f \big(2 \log \big(\frac{t+\sqrt{t^{2}+4}}{2}\big)\big)\big)}=\Big(\frac{t+\sqrt{t^{2}+4}}{2}\Big)^{2 x},
\end{align*}
where we used the fact that $\sum_{k=0}^{\infty}p_{k}(x)\frac{t^{k}}{k!}=e^{x \bar{f}(t)}$. This shows the assertion.\\
We would like to give another proof for this by using Theorem 5.1 (b). \\
Here we have to show that
\begin{align*}
\sum_{n=k}^{\infty}\frac{1}{k!}\langle{f(t)^{k}|x^{[n]}\rangle}\frac{t^{n}}{n!}=\frac{1}{k!}\big(f \big(2 \log \Big(\frac{t+\sqrt{t^{2}+4}}{2}\Big)\big)\big)^{k}.
\end{align*}
Equivalently, we need to verify that
\begin{align}
\sum_{n=k}^{\infty}\langle{f(t)^{k}|x^{[n]}\rangle}\frac{(e^{\frac{t}{2}}-e^{-\frac{t}{2}})^{n}}{n!}=(f(t))^{k}.\label{2E}
\end{align}
The left hand side of \eqref{2E} is equal to
\begin{align*}
\sum_{n=k}^{\infty}\sum_{l=k}^{n}T_{1}(n,l)\langle{f(t)^{k}|x^{l}\rangle}\frac{(e^{\frac{t}{2}}-e^{-\frac{t}{2}})^{n}}{n!}&=\sum_{l=k}^{\infty}\langle{f(t)^{k}|x^{l}\rangle}\sum_{n=l}^{\infty}T_{1}(n,l)\frac{(e^{\frac{t}{2}}-e^{-\frac{t}{2}})^{n}}{n!}\\
&=\sum_{l=0}^{\infty}\langle{f(t)^{k}|x^{l}\rangle}\frac{t^{l}}{l!}=(f(t))^{k}.
\end{align*}
\end{proof}
\end{theorem}

The following orthogonality and inverse relations can be easily shown from the definitions of central factorial numbers associated with $\bold{P}$.
\begin{theorem}
Let $\bold{P}=\left\{p_{n}(x)\right\}_{n=0}^{\infty}$ be a sequence of polynomials such that deg\,$p_{n}(x)=n, p_{0}(x)=1$. Then we have
\begin{align*}
&(a)\,\, \sum_{k=l}^{n}T_{1}(n,k;\bold{P})T_{2}(k,l;\bold{P})=\delta_{n,l},
\,\, \sum_{k=l}^{n}T_{2}(n,k;\bold{P})T_{1}(k,l;\bold{P})=\delta_{n,l},\,\,(0 \le l \le n). \\
&(b)\,\, a_{n}=\sum_{k=0}^{n}T_{2}(n,k;\bold{P})c_{k} \Longleftrightarrow c_{n}=\sum_{k=0}^{n}T_{1}(n,k;\bold{P})a_{k}.\\
&(c)\,\, a_{n}=\sum_{k=n}^{m}T_{2}(k,n;\bold{P})c_{k} \Longleftrightarrow c_{n}=\sum_{k=n}^{m}T_{1}(k,n;\bold{P})a_{k}.
\end{align*}
\end{theorem}

\begin{remark}
Let $p_{n}(x)$ be Sheffer for the pair $(g(t), f(t))$. Then we saw in Theorem 3.2 that the generating function of $T_{2}(n,k;\bold{P})$ is given by
\begin{align*}
\sum_{n=k}^{\infty}T_{2}(n,k;\bold{P})\frac{t^{n}}{n!}=\frac{1}{g(\bar{f}(t))}\frac{1}{k!}(e^{\frac{1}{2}\bar{f}(t)}-e^{-\frac{1}{2}\bar{f}(t)})^{k}.
\end{align*}
Let $f(t)$ be any delta series. After introducing the central logarithm associated to $f(t)$, it seems natural to introduce the {\it{central exponential associated to $f(t)$}} as 
\begin{align*}
EC_{f(t)}t=e^{\frac{1}{2}\bar{f}(t)}-e^{-\frac{1}{2}\bar{f}(t)}.
\end{align*}
When $p_{n}(x) \sim (g(t), f(t))$, $EC_{f(t)}t$ may be also denoted by $EC_{\bold{P}}t$ and called the central exponential associated to $\bold{P}$. Recalling that $x^{[n]} \sim (1,l(t)=e^{\frac{t}{2}}-e^{-\frac{t}{2}})$, we observe that $EC_{f(t)}t=l(\bar{f}(t))$. \\
As an example, consider the delta series $f(t)=\frac{1}{\lambda}(e^{\lambda t} -1)$. Here $\bar{f}(t)=\frac{1}{\lambda}\log(1+\lambda t)$. Thus we see that 
\begin{align}
&LC_{f(t)}t=\frac{1}{\lambda}\bigg(\bigg(\frac{t+\sqrt{t^{2}+4}}{2}\bigg)^{2 \lambda}-1 \bigg)=\log_{\lambda}\bigg(\frac{t+\sqrt{t^{2}+4}}{2}\bigg)^{2},\\
&EC_{f(t)}t=(1+\lambda t)^{\frac{1}{2 \lambda}}-(1+\lambda t)^{-\frac{1}{2 \lambda}}=e_{\lambda}^{\frac{1}{2}}(t)-e_{\lambda}^{-\frac{1}{2}}(t).
\end{align}
Here we note that, in many of the papers on degenerate special polynomials and numbers, $\frac{1}{\lambda}\big(t^{\lambda}-1 \big)$ is called the degenerate logarithm and denoted by $\log_{\lambda}t$, and $(1+\lambda t)^{\frac{1}{\lambda}}$ is called the degenerate exponential and denoted by $e_{\lambda}(t)$.
\end{remark}

\section{Examples on central factorial numbers of the first kind associated with sequences of polynomials}
(a) Let $\bold{P}=\left\{x^{n}\right\}$. Then $x^{n} \sim (1,t)$, and $T_{2}(n,k;\bold{P})=T_{2}(n,k)$.  As $x^{[n]}=\sum_{k=0}^{n}T_{1}(n,k;\bold{P})x^{k}$, we have (see [7,30,31])
\begin{align*}
T_{1}(n,k;\bold{P})=T_{1}(n,k), \quad \sum_{n=k}^{\infty}T_{1}(n,k)\frac{t^n}{n!}=\frac{1}{k!}\bigg(2\log\bigg(\frac{t+\sqrt{t^{2}+4}}{2}\bigg)\bigg)^{k},
\end{align*}
(see \eqref{1C}). Now, by Theorem 5.3, we get \\
\begin{align*}
&\sum_{k=l}^{n}T_{1}(n,k)T_{2}(k,l)=\delta_{n,l},
\,\, \sum_{k=l}^{n}T_{2}(n,k)T_{1}(k,l)=\delta_{n,l},\,\,(0 \le l \le n), \\
&a_{n}=\sum_{k=0}^{n}T_{2}(n,k)c_{k} \Longleftrightarrow c_{n}=\sum_{k=0}^{n}T_{1}(n,k)a_{k},\\
&a_{n}=\sum_{k=n}^{m}T_{2}(k,n)c_{k} \Longleftrightarrow c_{n}=\sum_{k=n}^{m}T_{1}(k,n)a_{k}.
\end{align*}

(b) Let $\bold{P}=\left\{(x)_{n,\lambda}\right\}$ be the sequence of generalized falling factorials. Then $(x)_{n,\lambda} \sim \big(1, f(t)=\frac{1}{\lambda}(e^{\lambda t}-1)\big)$, and $T_{2}(n,k;\bold{P})=T_{2,\lambda}(n,k)$. As $x^{[n]}=\sum_{k=0}^{n}T_{1}(n,k;\bold{P})(x)_{k,\lambda}$, we have (see [13,14,16,18,23-25])
\begin{align*}
T_{1}(n,k;\bold{P})=T_{1,\lambda}(n,k),\quad \sum_{n=k}^{\infty}T_{1,\lambda}(n,k)\frac{t^n}{n!}=\frac{1}{k!}\bigg(\log_{\lambda}\bigg(\frac{t+\sqrt{t^{2}+4}}{2}\bigg)^{2}\bigg)^{k},
\end{align*}
(see \eqref{3C}).  Now, by Theorem 5.3, we get \\
\begin{align*}
&\sum_{k=l}^{n}T_{1,\lambda}(n,k)T_{2,\lambda}(k,l)=\delta_{n,l},
\,\, \sum_{k=l}^{n}T_{2,\lambda}(n,k)T_{1,\lambda}(k,l)=\delta_{n,l},\,\,(0 \le l \le n), \\
&a_{n}=\sum_{k=0}^{n}T_{2,\lambda}(n,k)c_{k} \Longleftrightarrow c_{n}=\sum_{k=0}^{n}T_{1,\lambda}(n,k)a_{k},\\
&a_{n}=\sum_{k=n}^{m}T_{2,\lambda}(k,n)c_{k} \Longleftrightarrow c_{n}=\sum_{k=n}^{m}T_{1,\lambda}(k,n)a_{k}.
\end{align*}

(c) Let $\bold{P}=\left\{\langle{x \rangle}_{n}\right\}$ be the sequence of rising factorials. Then $\langle{x \rangle}_{n} \sim (1, 1-e^{-t})$, and $T_{2}(n,k;\bold{P})=\sum_{l=k}^{n}(-1)^{n-l}S_{1}(n,l)T_{2}(l,k)$ (see [7,15,30,31]). In addition, by Theorem 5.1, we have \\
\begin{align*}
T_1(n,k;\bold{P})&=\sum_{l=k}^{n}T_{1}(n,l)(-1)^{k}\big\langle{\frac{1}{k!}(e^{-t}-1)^{k}|x^{l}\big\rangle}\\
&=\sum_{l=k}^{n}T_{1}(n,l)(-1)^{k}\sum_{j=k}^{l}S_{2}(j,k)(-1)^{j}\frac{1}{j!}\langle{t^{j}|x^{l}\rangle}\\
&=\sum_{l=k}^{n}(-1)^{l-k}T_{1}(n,l)S_{2}(l,k).
\end{align*}
Thus $T_1(n,k;\bold{P})=\sum_{l=k}^{n}(-1)^{l-k}T_{1}(n,l)S_{2}(l,k)$.
Hence,  by Theorem 5.3,  we have the following identities given by:
\begin{align*}
&\sum_{k=l}^{n}\sum_{m=k}^{n}\sum_{j=l}^{k}T_{1}(n,m)(-1)^{m-k}S_{2}(m,k)(-1)^{k-j}S_{1}(k,j)T_{2}(j,l)=\delta_{n,l},\\
&\sum_{k=l}^{n}\sum_{j=k}^{n}\sum_{m=l}^{k}(-1)^{n-j}S_{1}(n,j)T_{2}(j,k)T_{1}(k,m)(-1)^{m-l}S_{2}(m,l)=\delta_{n,l}, \\
&a_{n}=\sum_{k=0}^{n}\sum_{l=k}^{n}(-1)^{n-l}S_{1}(n,l)T_{2}(l,k)c_{k} \Longleftrightarrow c_{n}=\sum_{k=0}^{n}\sum_{l=k}^{n}(-1)^{l-k}T_{1}(n,l)S_{2}(l,k)a_{k},\\
&a_{n}=\sum_{k=n}^{m}\sum_{l=n}^{k}(-1)^{k-l}S_{1}(k,l)T_{2}(l,n)c_{k} \Longleftrightarrow c_{n}=\sum_{k=n}^{m}\sum_{l=n}^{k}(-1)^{l-n}T_{1}(k,l)S_{2}(l,n)a_{k}. 
\end{align*}
(d) Let $\bold{P}=\left\{\langle{x \rangle}_{n,\lambda}\right\}$ be the sequence of generalized rising factorials. Then $\langle{x \rangle}_{n,\lambda} \sim \big(1, \frac{1}{\lambda}(1-e^{-\lambda t})\big)$, and $T_{2}(n,k;\bold{P})=\sum_{l=k}^{n}(-\lambda)^{n-l}S_{1}(n,l)T_{2}(l,k)$ (see [18]).
By Theorem 5.1, we have
\begin{align*}
T_{1}(n,k;\bold{P})&=\sum_{l=k}^{n}T_{1}(n,l)\frac{1}{(-\lambda)^k}\big\langle{\frac{1}{k!}(e^{-\lambda t}-1)^{k}|x^{l}\big\rangle}\\
&=\sum_{l=k}^{n}T_{1}(n,l)\frac{1}{(-\lambda)^k}\sum_{j=k}^{l}S_{2}(j,k)\frac{1}{j!}(-\lambda)^{j}\langle{t^{j}|x^{l}\rangle}\\
&=\sum_{l=k}^{n}(-\lambda)^{l-k}T_{1}(n,l)S_{2}(l,k).
\end{align*}
Thus $T_{1}(n,k;\bold{P})=\sum_{l=k}^{n}(-\lambda)^{l-k}T_{1}(n,l)S_{2}(l,k)$, (see [14]).  Hence,  by Theorem 5.3,  we get
\begin{align*}
&\sum_{k=l}^{n}\sum_{m=k}^{n}\sum_{j=l}^{k}T_{1}(n,m)(-\lambda)^{m-k}S_{2}(m,k)(-\lambda)^{k-j}S_{1}(k,j)T_{2}(j,l)=\delta_{n,l},\\
&\sum_{k=l}^{n}\sum_{j=k}^{n}\sum_{m=l}^{k}(-\lambda)^{n-j}S_{1}(n,j)T_{2}(j,k)T_{1}(k,m)(-\lambda)^{m-l}S_{2}(m,l)=\delta_{n,l},\\
&a_{n}=\sum_{k=0}^{n}\sum_{l=k}^{n}(-\lambda)^{n-l}S_{1}(n,l)T_{2}(l,k)c_{k} \Longleftrightarrow c_{n}=\sum_{k=0}^{n}\sum_{l=k}^{n}(-\lambda)^{l-k}T_{1}(n,l)S_{2}(l,k)a_{k},\\
&a_{n}=\sum_{k=0}^{n}\sum_{l=n}^{k}(-\lambda)^{k-l}S_{1}(k,l)T_{2}(l,n)c_{k} \Longleftrightarrow c_{n}=\sum_{k=0}^{n}\sum_{l=n}^{k}(-\lambda)^{l-n}T_{1}(k,l)S_{2}(l,n)a_{k}. 
\end{align*}

(e) Let $\bold{P_{1}}=\left\{TLB_{n,1}(x)\right\}_{n=0}^{\infty}$ be the sequence of central factorial Lah Bell polynomials of the first kind. Then we saw that $T_{2}(n,k;\bold{P_{1}})=L_{2}^{(c)}(n,k)$, where $L_{2}^{(c)}(n,k)$ are the central Lah numbers of the second kind (see [19]). We note from \eqref{6D} that
\begin{align*}
\sum_{n=0}^{\infty}x^{[n]}\frac{t^{n}}{n!}&=e^{x\bar{f}(t)}=\sum_{k=0}^{\infty}TLB_{k,1}(x)\frac{1}{k!}(\bar{\alpha}(t))^{k}\\
&=\sum_{k=0}^{\infty}TLB_{k,1}(x)\sum_{n=k}^{\infty}L_{1}^{(c)}(n,k)\frac{t^{n}}{n!}\\
&=\sum_{n=0}^{\infty}\sum_{k=0}^{n}L_{1}^{(c)}(n,k)TLB_{k,1}(x)\frac{t^{n}}{n!}.
\end{align*}
Thus we have shown that
\begin{align*}
x^{[n]}=\sum_{k=0}^{n}L_{1}^{(c)}(n,k)TLB_{k,1}(x),
\end{align*}
which verifies that $T_{1}(n,k;\bold{P_{1}})=L_{1}^{(c)}(n,k)$.
Now, by Theorem 5.3, we get \\
\begin{align*}
&\sum_{k=l}^{n}L_{1}^{(c)}(n,k)L_{2}^{(c)}(k,l)=\delta_{n,l},
\,\, \sum_{k=l}^{n}L_{2}^{(c)}(n,k)L_{1}^{(c)}(k,l)=\delta_{n,l},\,\,(0 \le l \le n), \\
&a_{n}=\sum_{k=0}^{n}L_{2}^{(c)}(n,k)c_{k} \Longleftrightarrow c_{n}=\sum_{k=0}^{n}L_{1}^{(c)}(n,k)a_{k},\\
&a_{n}=\sum_{k=n}^{m}L_{2}^{(c)}(k,n)c_{k} \Longleftrightarrow c_{n}=\sum_{k=n}^{m}L_{1}^{(c)}(k,n)a_{k}.
\end{align*}
Now, we let $\bold{P_{2}}=\left\{TLB_{n,2}(x)\right\}_{n=0}^{\infty}$ (see \eqref{8D}). Then we saw that $T_{2}(n,k;\bold{P_{2}}) = \sum_{l=k}^{n}TL_{2}(n,l)T_{2}(l,k)$, and $TLB_{n,2}(x) \sim (1, \bar{f}(\alpha(t))$.
Recalling from \eqref{5D} that $\frac{1}{k!}\big(\bar{f}(\alpha(t))\big)^{k}=\sum_{n=k}^{\infty}TL_{1}(n,k)\frac{t^{n}}{n!}$, from Theorem 5.1 we have
\begin{align*}
T_{1}(n,k;\bold{P_{2}})&=\sum_{l=k}^{n}T_{1}(n,l)\big\langle{\frac{1}{k!}\big(\bar{f}(\alpha(t))\big)^{k} \big|x^{l}\big\rangle}\\
&=\sum_{l=k}^{n}T_{1}(n,l)\big\langle{\sum_{m=k}^{\infty}TL_{1}(m,k)\frac{t^{m}}{m!} \big|x^{l}\big \rangle}\\
&=\sum_{l=k}^{n}T_{1}(n,l)\sum_{m=k}^{l}TL_{1}(m,k)\frac{1}{m!}\langle{t^{m}|x^{l}\rangle}\\
&=\sum_{l=k}^{n}T_{1}(n,l)TL_{1}(l,k).
\end{align*}
Thus we have shown that $T_{1}(n,k;\bold{P_{2}})=\sum_{l=k}^{n}T_{1}(n,l)TL_{1}(l,k)$.
Now, by Theorem 5.3 we obtain
\begin{align*}
&\sum_{k=l}^{n}\sum_{m=k}^{n}\sum_{j=l}^{k}T_{1}(n,m)TL_{1}(m,k)TL_{2}(k,j)T_{2}(j,l)=\delta_{n,l}, \\
&\sum_{k=l}^{n}\sum_{j=k}^{n}\sum_{m=l}^{k}TL_{2}(n,j)T_{2}(j,k)T_{1}(k,m)TL_{1}(m,l)=\delta_{n,l},\,\,(0 \le l \le n), \\
&a_{n}=\sum_{k=0}^{n} \sum_{l=k}^{n}TL_{2}(n,l)T_{2}(l,k)c_{k} \Longleftrightarrow c_{n}=\sum_{k=0}^{n}\sum_{l=k}^{n}T_{1}(n,l)TL_{1}(l,k)a_{k},\\
&a_{n}=\sum_{k=n}^{m}\sum_{l=n}^{k}TL_{2}(k,l)T_{2}(l,n)c_{k} \Longleftrightarrow c_{n}=\sum_{k=n}^{m}\sum_{l=n}^{k}T_{1}(k,l)TL_{1}(l,n)a_{k}.
\end{align*}
(f) Let $\bold{P}=\left\{\mathrm{Bel}_{n}^{(c)}(x)\right\}$ be the sequence of central Bell polynomials. Then $\mathrm{Bel}_{n}^{(c)}(x) \sim \Big(1, f(t)=2\log\big(\frac{t+\sqrt{t^{2}+4}}{2}\big)\Big)$, and $T_{2}(n,k;\bold{P})=\sum_{l=k}^{n}T_{2}(n,l)T_{2}(l,k)$ (see [22]).
Now, by Theorem 5.1 and \eqref{1C}, we get
\begin{align*}
T_{1}(n,k;\bold{P})&=\sum_{l=k}^{n}T_{1}(n,l)\big\langle{\frac{1}{k!}\Big(2 \log\Big(\frac{t+\sqrt{t^{2}+4}}{2}\Big)\Big)^{k}|x^{l}\big\rangle}\\
&=\sum_{l=k}^{n}T_{1}(n,l)\sum_{j=k}^{l}T_{1}(j,k)\frac{1}{j!}\langle{t^{j}|x^{l}\rangle}\\
&=\sum_{l=k}^{n}T_{1}(n,l)T_{1}(l,k).
\end{align*}
Therefore $T_{1}(n,k;\bold{P})=\sum_{l=k}^{n}T_{1}(n,l)T_{1}(l,k)$ (see [15]).  From Theorem 5.3,  we have:
\begin{align*}
&\sum_{k=l}^{n}\sum_{m=k}^{n}\sum_{j=l}^{k}T_{1}(n,m)T_{1}(m,k)T_{2}(k,j)T_{2}(j,l)=\delta_{n,l},\\
&\sum_{k=l}^{n}\sum_{j=k}^{n}\sum_{m=l}^{k}T_{2}(n,j)T_{2}(j,k)T_{1}(k,m)T_{1}(m,l)=\delta_{n,l},\\ 
&a_{n}=\sum_{k=0}^{n}\sum_{l=k}^{n}T_{2}(n,l)T_{2}(l,k)c_{k} \Longleftrightarrow c_{n}=\sum_{k=0}^{n}\sum_{l=k}^{n}T_{1}(n,l)T_{1}(l,k)a_{k},\\
&a_{n}=\sum_{k=n}^{m}\sum_{l=n}^{k}T_{2}(k,l)T_{2}(l,n)c_{k} \Longleftrightarrow c_{n}=\sum_{k=n}^{m}\sum_{l=n}^{k}T_{1}(k,l)T_{1}(l,n)a_{k}.
\end{align*}

(g) Let $\bold{P}=\left\{\mathrm{Bel}_{n,\lambda}^{(c)}(x)\right\}$ be the sequence of degenerate central Bell polynomials. Then $\mathrm{Bel}_{n,\lambda}^{(c)}(x) \sim \Big(1, f(t)=\log_{\lambda}\Big(\frac{t+\sqrt{t^{2}+4}}{2}\Big)^{2}\Big),$ and $T_{2}(n,k;\bold{P})=\sum_{l=k}^{n}T_{2,\lambda}(n,l)T_{2}(l,k)$ (see [23]).
Now, by Theorem 5.1 and \eqref{3C}, we have
\begin{align*}
T_{1}(n,k;\bold{P})&=\sum_{l=k}^{n}T_{1}(n,l)\big\langle{\frac{1}{k!}\Big(\log_{\lambda}\Big(\frac{t+\sqrt{t^{2}+4}}{2}\Big)^{2}\Big)^{k}|x^{l}\big\rangle} \\
&=\sum_{l=k}^{n}T_{1}(n,l)\sum_{j=k}^{l}T_{1,\lambda}(j,k)\frac{1}{j!}\langle{t^{j}|x^{l}\rangle}\\
&=\sum_{l=k}^{n}T_{1}(n,l)T_{1,\lambda}(l,k).
\end{align*}
Thus we have $T_{1}(n,k;\bold{P})=\sum_{l=k}^{n}T_{1}(n,l)T_{1,\lambda}(l,k)$ (see [16]).
From Theorem 5.3, we have:
\begin{align*}
&\sum_{k=l}^{n}\sum_{m=k}^{n}\sum_{j=l}^{k}T_{1}(n,m)T_{1,\lambda}(m,k)T_{2,\lambda}(k,j)T_{2}(j,l)=\delta_{n,l},\\
&\sum_{k=l}^{n}\sum_{j=k}^{n}\sum_{m=l}^{k}T_{2,\lambda}(n,j)T_{2}(j,k)T_{1}(k,m)T_{1,\lambda}(m,l)=\delta_{n,l},\\
&a_{n}=\sum_{k=0}^{n}\sum_{l=k}^{n}T_{2,\lambda}(n,l)T_{2}(l,k)c_{k} \Longleftrightarrow c_{n}=\sum_{k=0}^{n}\sum_{l=k}^{n}T_{1}(n,l)T_{1,\lambda}(l,k)a_{k},\\
&a_{n}=\sum_{k=n}^{m}\sum_{l=n}^{k}T_{2,\lambda}(k,l)T_{2}(l,n)c_{k} \Longleftrightarrow c_{n}=\sum_{k=n}^{m}\sum_{l=n}^{k}T_{1}(k,l)T_{1,\lambda}(l,n)a_{k}.
\end{align*}

(h) Let $x^{[n,\lambda]}$ be the sequence of polynomials defined by $x^{[0,\lambda]}=1,\quad x^{[n,\lambda]}=x\big(x+(\frac{1}{2}n-1)\lambda\big)_{n-1,\lambda}, \,\,(n \ge 1)$. Then $x^{[n,\lambda]} \sim \big(1, f(t)=\frac{1}{\lambda}(e^{\frac{\lambda t}{2}}-e^{-\frac{\lambda t}{2}})\big)$, and $T_{2}(n,k;\bold{P})=\sum_{l=k}^{n}R_{1,\lambda}(n,l)T_{2}(l,k)$,
where $R_{1,\lambda}(n,k)$ are defined by $x^{[n,\lambda]}=\sum_{k=0}^{n}R_{1,\lambda}(n,k)x^{k}$. Here we recall that
\begin{align*}
 \sum_{n=k}^{\infty}R_{2,\lambda}(n,k)\frac{t^{n}}{n!}=\frac{1}{k!}\Big(\frac{1}{\lambda}\big(e^{\frac{\lambda t}{2}}-e^{-\frac{\lambda t}{2}}\big)\Big)^{k}.
\end{align*}
Let $\bold{P}=\left\{x^{[n,\lambda]}\right\}$. Now, by Theorem 5.1, we have
\begin{align*}
T_{1}(n,k;\bold{P})&=\sum_{l=k}^{n}T_{1}(n,l)\big\langle{\frac{1}{k!}\big(\frac{1}{\lambda}\big(e^{\frac{\lambda t}{2}}-e^{-\frac{\lambda t}{2}}\big)\big)^{k}|x^{l}\big\rangle}\\
&=\sum_{l=k}^{n}T_{1}(n,l)\sum_{j=k}^{l}R_{2,\lambda}(j,k)\frac{1}{j!}\langle{t^{j}|x^{l}\rangle}\\
&=\sum_{l=k}^{n}T_{1}(n,l)R_{2,\lambda}(l,k).
\end{align*}
Hence $T_{1}(n,k;\bold{P})=\sum_{l=k}^{n}T_{1}(n,l)R_{2,\lambda}(l,k)$. From Theorem 5.3, we get
\begin{align*}
&\sum_{k=l}^{n}\sum_{m=k}^{n}\sum_{j=l}^{k}T_{1}(n,m)R_{2,\lambda}(m,k)R_{1,\lambda}(k,j)T_{2}(j,l)=\delta_{n,l}, \\
&\sum_{k=l}^{n}\sum_{j=k}^{n}\sum_{m=l}^{k}R_{1,\lambda}(n,j)T_{2}(j,k)T_{1}(k,m)R_{2,\lambda}(m,l)=\delta_{n,l},\\
&a_{n}=\sum_{k=0}^{n}\sum_{l=k}^{n}R_{1,\lambda}(n,l)T_{2}(l,k)c_{k} \Longleftrightarrow c_{n}=\sum_{k=0}^{n}\sum_{l=k}^{n}T_{1}(n,l)R_{2,\lambda}(l,k)a_{k},\\
&a_{n}=\sum_{k=n}^{m}\sum_{l=n}^{k}R_{1,\lambda}(k,l)T_{2}(l,n)c_{k} \Longleftrightarrow c_{n}=\sum_{k=n}^{m}\sum_{l=n}^{k}T_{1}(k,l)R_{2,\lambda}(l,n)a_{k}.
\end{align*}

(i) Let $\bold{P}=\left\{B_{n}^{L}(x)\right\}$ be the sequence of Lah-Bell polynomials. Then $B_{n}^{L}(x) \sim (1,\frac{t}{1+t})$, and $T_{2}(n,k;\bold{P})=\sum_{l=k}^{n}L(n,l)T_{2}(l,k)$ (see [7,15,30,31]).
By Theorem 5.1, we see that
\begin{align*}
T_{1}(n,k;\bold{P})&=\sum_{l=k}^{n}T_{1}(n,l)(-1)^{k}\big\langle{\frac{1}{k!}\Big(\frac{-t}{1-(-t)}\Big)^{k}|x^{l}\big\rangle}\\
&=\sum_{l=k}^{n}T_{1}(n,l)(-1)^{k}\sum_{j=k}^{l}L(j,k)(-1)^{j}\frac{1}{j!}\langle{t^{j}|x^{l}\rangle}\\
&=\sum_{l=k}^{n}(-1)^{l-k}T_{1}(n,l)L(l,k).
\end{align*}
Therefore we have $T_{1}(n,k;\bold{P})=\sum_{l=k}^{n}(-1)^{l-k}T_{1}(n,l)L(l,k)$ (see [11]).  Now,  from Theorem 5.3, we obtain 
\begin{align*}
&\sum_{k=l}^{n}\sum_{m=k}^{n}\sum_{j=l}^{k}(-1)^{m-k}T_{1}(n,m)L(m,k)L(k,j)T_2(j,l)=\delta_{n,l}, \\
&\sum_{k=l}^{n}\sum_{j=k}^{n}\sum_{m=l}^{k}(-1)^{m-l}L(n,j)T_{2}(j,k)T_{1}(k,m)L(m,l)=\delta_{n,l},\\
&a_{n}=\sum_{k=0}^{n}\sum_{l=k}^{n}L(n,l)T_{2}(l,k)c_{k} \Longleftrightarrow c_{n}=\sum_{k=0}^{n}\sum_{l=k}^{n}(-1)^{l-k}T_{1}(n,l)L(l,k)a_{k},\\
&a_{n}=\sum_{k=n}^{m}\sum_{l=n}^{k}L(k,l)T_{2}(l,n)c_{k} \Longleftrightarrow c_{n}=\sum_{k=n}^{l}\sum_{l=n}^{k}(-1)^{l-n}T_{1}(k,l)L(l,n)a_{k}.
\end{align*}
(j) Let $\bold{P}=\left\{B_{n,\lambda}^{L}(x)\right\}$ be the sequence of degenerate Lah-Bell polynomials. Then $B_{n,\lambda}^{L}(x) \sim (1,\frac{e^{\lambda t}-1}{\lambda+e^{\lambda t}-1})$, and $T_{2}(n,k;\bold{P})=\sum_{l=k}^{n}\sum_{m=l}^{n}\lambda^{m-l}L(n,m)S_{1}(m,l)T_{2}(l,k).$ Now, by Theorem 5.1, we get
\begin{align*}
T_{1}(n,k;\bold{P})&=\sum_{l=k}^{n}T_{1}(n,l)(-1)^{k}\Big\langle{\frac{1}{k!}\Big(\frac{-\frac{1}{\lambda}(e^{\lambda t}-1)}{1-(-\frac{1}{\lambda}(e^{\lambda t}-1))}\Big)^{k}|x^{l}\Big\rangle}\\
&=\sum_{l=k}^{n}T_{1}(n,l)(-1)^{k}\sum_{m=k}^{l}L(m,k)\big(-\frac{1}{\lambda}\big)^{m}\big\langle{\frac{1}{m!}(e^{\lambda t}-1)^{m}|x^{l}\big\rangle}\\
&=\sum_{l=k}^{n}T_{1}(n,l)(-1)^{k}\sum_{m=k}^{l}L(m,k)\big(-\frac{1}{\lambda}\big)^{m}
\sum_{j=m}^{l}                                                                                                                                                                                                                                                                                                                                                                                                                                                                                                                                                                  S_{2}(j,m)\lambda^{j}\frac{1}{j!}\langle{t^{j}|x^{l}\rangle}\\
&=\sum_{l=k}^{n}\sum_{m=k}^{l}(-1)^{m-k}\lambda^{l-m}T_{1}(n,l)S_{2}(l,m)L(m,k).
\end{align*}
Thus we have $T_{1}(n,k;\bold{P})=\sum_{l=k}^{n}\sum_{m=k}^{l}(-1)^{m-k}\lambda^{l-m}T_{1}(n,l)S_{2}(l,m)L(m,k)$, (see [14]).
Now, by Theorem 5.3, we obtain 
\begin{align*}
&a_{n}=\sum_{k=0}^{n}\sum_{l=k}^{n}\sum_{m=l}^{n}\lambda^{m-l}L(n,m)S_{1}(m,l)T_{2}(l,k)c_{k} \\ 
&\Longleftrightarrow c_{n}=\sum_{k=0}^{n}\sum_{l=k}^{n}\sum_{m=k}^{l}(-1)^{m-k}\lambda^{l-m}T_{1}(n,l)S_{2}(l,m)L(m,k)a_{k},\\
&a_{n}=\sum_{k=n}^{m}\sum_{l=n}^{k}\sum_{m=l}^{k}\lambda^{m-l}L(k,m)S_{1}(m,l)T_{2}(l,n)c_{k}\\
&\Longleftrightarrow c_{n}=\sum_{k=n}^{m}\sum_{l=n}^{k}\sum_{m=n}^{l}(-1)^{m-n}\lambda^{l-m}T_{1}(k,l)S_{2}(l,m)L(m,n)a_{k}.
\end{align*}

(k) Let $\bold{P}=\left\{\mathrm{Bel}_{n}(x)\right\}$ be the sequence of Bell polynomials. Then $\mathrm{Bel}_{n}(x) \sim (1, \log (1+t))$, and $T_{2}(n,k;\bold{P})=\sum_{l=k}^{n}S_2(n,l)T_{2}(l,k)$ 
(see [7,30,31]). By Theorem 5.1, we get
\begin{align*}
T_{1}(n,k;\bold{P})&=\sum_{l=k}^{n}T_{1}(n,l)\big\langle{\frac{1}{k!}(\log (1+t))^{k}|x^{l}\big\rangle}\\
&=\sum_{l=k}^{n}T_{1}(n,l)\sum_{j=k}^{l}S_{1}(j,k)\frac{1}{j!}\langle{t^{j}|x^{l}\rangle}\\
&=\sum_{l=k}^{n}T_{1}(n,l)S_{1}(l,k).
\end{align*}
Thus we see that $T_{1}(n,k;\bold{P})=\sum_{l=k}^{n}T_{1}(n,l)S_{1}(l,k)$ (see [4,22,23]).  Now,  by Theorem 5.3,  we have
\begin{align*}
&\sum_{k=l}^{n}\sum_{m=k}^{n}\sum_{j=l}^{k}T_{1}(n,m)S_{1}(m,k)S_{2}(k,j)T_{2}(j,l)=\delta_{n,l}, \\
&\sum_{k=l}^{n}\sum_{j=k}^{n}\sum_{m=l}^{k}S_{2}(n,j)T_{2}(j,k)T_{1}(k,m)S_{1}(m,l)=\delta_{n,l},\\
&a_{n}=\sum_{k=0}^{n}\sum_{l=k}^{n}S_2(n,l)T_{2}(l,k)c_{k} \Longleftrightarrow c_{n}=\sum_{k=0}^{n}\sum_{l=k}^{n}T_{1}(n,l)S_{1}(l,k)a_{k},\\
&a_{n}=\sum_{k=n}^{m}\sum_{l=n}^{k}S_2(k,l)T_{2}(l,n)c_{k} \Longleftrightarrow c_{n}=\sum_{k=n}^{m}\sum_{l=n}^{k}T_{1}(k,l)S_{1}(l,n)a_{k}.
\end{align*}

(l) Let $\bold{P_{1}}=\left\{\mathrm{Bel}_{n,\lambda}(x)\right\}$ be the sequence of partially degenerate Bell polynomials. Then $\mathrm{Bel}_{n,\lambda}(x) \sim (1, \log_{\lambda}(1+t))$, and 
$T_{2}(n,k;\bold{P_{1}})=\sum_{l=k}^{n}S_{2,\lambda}(n,l)T_{2}(l,k)$ (see [24]). By Theorem 5.1, we have
\begin{align*}
T_{1}(n,k;\bold{P_{1}})&=\sum_{l=k}^{n}T_{1}(n,l)\big\langle{\frac{1}{k!}(\log_{\lambda} (1+t))^{k}|x^{l}\big\rangle}\\
&=\sum_{l=k}^{n}T_{1}(n,l)\sum_{j=k}^{l}S_{1,\lambda}(j,k)\frac{1}{j!}\langle{t^{j}|x^{l}\rangle}\\
&=\sum_{l=k}^{n}T_{1}(n,l)S_{1,\lambda}(l,k).
\end{align*}
Hence we have $T_{1}(n,k;\bold{P_{1}})=\sum_{l=k}^{n}T_{1}(n,l)S_{1,\lambda}(l,k)$. Now, by Theorem 5.3,  we obtain 
\begin{align*}
&\sum_{k=l}^{n}\sum_{m=k}^{n}\sum_{j=l}^{k}T_{1}(n,m)S_{1,\lambda}(m,k)S_{2,\lambda}(k,j)T_{2}(j,l)=\delta_{n,l}, \\
&\sum_{k=l}^{n}\sum_{j=k}^{n}\sum_{m=l}^{k}S_{2,\lambda}(n,j)T_{2}(j,k)T_{1}(k,m)S_{1,\lambda}(m,l)=\delta_{n,l},\\
&a_{n}=\sum_{k=0}^{n}\sum_{l=k}^{n}S_{2,\lambda}(n,l)T_{2}(l,k)c_{k} \Longleftrightarrow c_{n}=\sum_{k=0}^{n}\sum_{l=k}^{n}T_{1}(n,l)S_{1,\lambda}(l,k)a_{k},\\
&a_{n}=\sum_{k=n}^{m}\sum_{l=n}^{k}S_{2,\lambda}(k,l)T_{2}(l,n)c_{k} \Longleftrightarrow c_{n}=\sum_{k=n}^{m}\sum_{l=n}^{k}T_{1}(k,l)S_{1,\lambda}(l,n)a_{k}.
\end{align*}

Let $\bold{P_{2}}=\left\{\phi_{n,\lambda}(x)\right\}$ be the sequence of fully degenerate Bell polynomials. Then $\phi_{n,\lambda}(x)$ $\sim \big(1, \log_{\lambda}(1+\frac{1}{\lambda}(e^{\lambda t}-1))\big)$, and $T_{2}(n,k;\bold{P_{2}})=\sum_{l=k}^{n}\sum_{m=l}^{n}\lambda^{m-l}S_{2,\lambda}(n,m)S_{1}(m,l)T_{2}(l,k)$ (see [25]). By Theorem 5.1, we obtain
\begin{align*}
T_{1}(n,k;\bold{P_{2}})&=\sum_{l=k}^{n}T_{1}(n,l)\big\langle{\frac{1}{k!}(\log_{\lambda}(1+\frac{1}{\lambda}(e^{\lambda t}-1)))^{k}|x^{l}\rangle}\\
&=\sum_{l=k}^{n}T_{1}(n,l)\sum_{m=k}^{l}S_{1,\lambda}(m,k)\frac{1}{\lambda^{m}}\langle{\frac{1}{m!}(e^{\lambda t}-1)^{m}|x^{l}\rangle}\\
&=\sum_{l=k}^{n}T_{1}(n,l)\sum_{m=k}^{l}S_{1,\lambda}(m,k)\frac{1}{\lambda^{m}}\sum_{j=m}^{l}S_{1}(j,m)\lambda^{j}\frac{1}{j!}\langle{t^{j}|x^{l}\rangle}\\
&=\sum_{l=k}^{n}\sum_{m=k}^{l}\lambda^{l-m}T_{1}(n,l)S_{1}(l,m)S_{1,\lambda}(m,k).
\end{align*}
Therefore we get $T_{1}(n,k;\bold{P_{2}})=\sum_{l=k}^{n}\sum_{m=k}^{l}\lambda^{l-m}T_{1}(n,l)S_{1}(l,m)S_{1,\lambda}(m,k)$ (see [18]).
Now, by Theorem 5.3, we obtain 
\begin{align*}
&a_{n}=\sum_{k=0}^{n}\sum_{l=k}^{n}\sum_{m=l}^{n}\lambda^{m-l}S_{2,\lambda}(n,m)S_{1}(m,l)T_{2}(l,k)c_{k} \\
&\Longleftrightarrow c_{n}=\sum_{k=0}^{n}\sum_{l=k}^{n}\sum_{m=k}^{l}\lambda^{l-m}T_{1}(n,l)S_{1}(l,m)S_{1,\lambda}(m,k)a_{k},\\
&a_{n}=\sum_{k=n}^{m}\sum_{l=n}^{k}\sum_{m=l}^{k}\lambda^{m-l}S_{2,\lambda}(k,m)S_{1}(m,l)T_{2}(l,n)c_{k} \\
&\Longleftrightarrow c_{n}=\sum_{k=n}^{m}\sum_{l=n}^{k}\sum_{m=n}^{l}\lambda^{l-m}T_{1}(k,l)S_{1}(l,m)S_{1,\lambda}(m,n)a_{k}.
\end{align*}

(m) Let $\bold{P}=\left\{M_{n}(x)\right\}$ be the sequence of Mittag-Leffler polynomials. Then $M_{n}(x) \sim (1,f(t)=\frac{e^{t}-1}{e^{t}+1})$, and $T_{2}(n,k;\bold{P})=\sum_{l=k}^{n}\sum_{m=l}^{n}2^{m}L(n,m)S_{1}(m,l)T_{2}(l,k)$ (see [31]). Then, from Theorem 5.1, we have
\begin{align*}
T_{1}(n,k;\bold{P})&=\sum_{l=k}^{n}T_{1}(n,l)\big\langle{\frac{1}{k!}\Big(\frac{e^{t}-1}{e^{t}+1}\Big)^{k}|x^{l}\big\rangle}\\
&=\sum_{l=k}^{n}T_{1}(n,l)(-1)^{k}\Big\langle{\frac{1}{k!}\Big(\frac{-\frac{1}{2}(e^{t}-1)}{1-(-\frac{1}{2}(e^{t}-1))}\Big)^{k}\Big|x^{l}\Big\rangle}\\
&=\sum_{l=k}^{n}T_{1}(n,l)(-1)^{k}\sum_{m=k}^{l}L(m,k)(-\frac{1}{2})^{m}\big\langle{\frac{1}{m!}(e^{t}-1)^{m}|x^{l}\big\rangle}\\
&=\sum_{l=k}^{n}T_{1}(n,l)(-1)^{k}\sum_{m=k}^{l}L(m,k)(-\frac{1}{2})^{m}\sum_{j=m}^{l}S_{2}(j,m)\frac{1}{j!}\langle{t^{j}|x^{l}\rangle}\\
&=\sum_{l=k}^{n}\sum_{m=k}^{l}(-1)^{m-k}2^{-m}T_{1}(n,l)S_{2}(l,m)L(m,k).
\end{align*}
Hence we have $T_{1}(n,k;\bold{P})=\sum_{l=k}^{n}\sum_{m=k}^{l}(-1)^{m-k}2^{-m}T_{1}(n,l)S_{2}(l,m)L(m,k)$. By Theorem 5.3, we get
\begin{align*}
&a_{n}=\sum_{k=0}^{n}\sum_{l=k}^{n}\sum_{m=l}^{n}2^{m}L(n,m)S_{1}(m,l)T_{2}(l,k)c_{k} \\ &\Longleftrightarrow c_{n}=\sum_{k=0}^{n}\sum_{l=k}^{n}\sum_{m=k}^{l}(-1)^{m-k}2^{-m}T_{1}(n,l)S_{2}(l,m)L(m,k)a_{k},\\
&a_{n}=\sum_{k=n}^{m}\sum_{l=n}^{k}\sum_{m=l}^{k}2^{m}L(k,m)S_{1}(m,l)T_{2}(l,n)c_{k} \\ &\Longleftrightarrow c_{n}=\sum_{k=n}^{m}\sum_{l=n}^{k}\sum_{m=n}^{l}(-1)^{m-n}2^{-m}T_{1}(k,l)S_{2}(l,m)L(m,n)a_{k}.
\end{align*}

(n) Let $\bold{P}=\left\{L_{n}(x)\right\}$ be the sequence of Laguerre polynomials of order -1.
Then $L_{n}(x) \sim (1, f(t)=\frac{t}{t-1})$, and $T_{2}(n,k;\bold{P})=\sum_{l=k}^{n}(-1)^{l}L(n,l)T_{2}(l,k)$ (see [31]). By Theorem 5.1, we see that
\begin{align*}
T_{1}(n,k;\bold{P})&=\sum_{l=k}^{n}T_{1}(n,l)\big\langle{\frac{1}{k!}\Big(\frac{t}{t-1}\Big)^{k}|x^{l}\big\rangle}\\
&=\sum_{l=k}^{n}T_{1}(n,l)(-1)^{k}\big\langle{\frac{1}{k!}\Big(\frac{t}{1-t}\Big)^{k}|x^{l}\big\rangle}\\
&=\sum_{l=k}^{n}T_{1}(n,l)(-1)^{k}\sum_{j=k}^{l}L(j,k)\frac{1}{j!}\langle{t^{j}|x^{l}\rangle}\\
&=(-1)^{k}\sum_{l=k}^{n}T_{1}(n,l)L(l,k).
\end{align*}
Thus we have $T_{1}(n,k;\bold{P})=(-1)^{k}\sum_{l=k}^{n}T_{1}(n,l)L(l,k)$ (see [23]).  Now, from Theorem 5.3, we get
\begin{align*}
&\sum_{k=l}^{n}\sum_{m=k}^{n}\sum_{j=l}^{k}(-1)^{k-j}T_{1}(n,m)L(m,k)L(k,j)T_{2}(j,l)=\delta_{n,l},\\ 
&\sum_{k=l}^{n}\sum_{j=k}^{n}\sum_{m=l}^{k}(-1)^{j-l}L(n,j)T_{2}(j,k)T_{1}(k,m)L(m,l)=\delta_{n,l}, \\
&a_{n}=\sum_{k=0}^{n}\sum_{l=k}^{n}(-1)^{l}L(n,l)T_{2}(l,k)c_{k} \Longleftrightarrow c_{n}=\sum_{k=0}^{n}(-1)^{k}\sum_{l=k}^{n}T_{1}(n,l)L(l,k)a_{k},\\
&a_{n}=\sum_{k=n}^{m}\sum_{l=n}^{k}(-1)^{l}L(k,l)T_{2}(l,n)c_{k} \Longleftrightarrow c_{n}=\sum_{k=n}^{m}(-1)^{n}\sum_{l=n}^{k}T_{1}(k,l)L(l,n)a_{k}.
\end{align*}

(o) Let $\bold{P}=\left\{B_{n}(x)\right\}$ be the sequence of Bernoulli polynomials. Then $B_{n}(x) \sim \big(\frac{e^{t}-1}{t},t\big)$, and $T_{2}(n,k;\bold{P})=\sum_{l=k}^{n}T_{2}(l,k)\binom{n}{l}B_{n-l}$ (see [31]). By using \eqref{7B} and Theorem 5.1, we see that 
\begin{align*}
T_{1}(n,k;\bold{P})&=\frac{1}{k!}\sum_{l=k}^{n}T_{1}(n,l)\Big\langle{t^{k}\Big|\frac{e^{t}-1}{t}x^{l}\Big\rangle}\\
&=\frac{1}{k!}\sum_{l=k}^{n}T_{1}(n,l)\Big\langle{t^{k}\Big|\int_{x}^{x+1}u^{l}du\Big\rangle}\\
&=\frac{1}{k!}\sum_{l=k}^{n}T_{1}(n,l)\frac{1}{l+1}(l+1)_{k}\Big\langle{ 1\Big|(x+1)^{l+1-k}-x^{l+1-k}\Big\rangle}\\
&=\sum_{l=k}^{n}\frac{1}{l+1}\binom{l+1}{k}T_{1}(n,l).
\end{align*}
Thus we have $T_{1}(n,k;\bold{P})=\sum_{l=k}^{n}\frac{1}{l+1}\binom{l+1}{k}T_{1}(n,l)$.
Now, by Theorem 5.3, we obtain
\begin{align*}
&\sum_{k=l}^{n}\sum_{m=k}^{n}\sum_{j=l}^{k}\frac{1}{m+1}\binom{k}{j}\binom{m+1}{k}T_{1}(n,m)T_{2}(j,l)B_{k-j}=\delta_{n,l},\\
&\sum_{k=l}^{n}\sum_{j=k}^{n}\sum_{m=l}^{k}\frac{1}{m+1}\binom{n}{j}\binom{m+1}{l}T_{2}(j,k)T_{1}(k,m)B_{n-j}=\delta_{n,l},\\
&a_{n}=\sum_{k=0}^{n}\sum_{l=k}^{n}T_{2}(l,k)\binom{n}{l}B_{n-l}c_{k} \Longleftrightarrow c_{n}=\sum_{k=0}^{n}\sum_{l=k}^{n}\frac{1}{l+1}\binom{l+1}{k}T_{1}(n,l)a_{k},\\
&a_{n}=\sum_{k=0}^{n}\sum_{l=n}^{k}T_{2}(l,n)\binom{k}{l}B_{k-l}c_{k} \Longleftrightarrow c_{n}=\sum_{k=0}^{n}\sum_{l=n}^{k}\frac{1}{l+1}\binom{l+1}{n}T_{1}(k,l)a_{k}.
\end{align*}

(p) Let $\bold{P}=\left\{E_{n}(x)\right\}$ be the sequence of Euler polynomials. Then $E_{n}(x) \sim \big(\frac{e^{t}+1}{2},t\big)$, and $T_{2}(n,k;\bold{P})=\sum_{l=k}^{n}\binom{n}{l}T_{2}(l,k)E_{n-l}$ (see [31]). By Theorem 5.1, we note that 
\begin{align*}
T_{1}(n,k;\bold{P})&=\frac{1}{2 k!}\sum_{l=k}^{n}T_{1}(n,l)\big\langle{e^{t}+1\big| t^{k}x^{l}\big\rangle}
=\frac{1}{2 k!}\sum_{l=k}^{n}T_{1}(n,l)(l)_{k}\big\langle{e^{t}+1\big| x^{l-k}
\rangle}\\
&=\frac{1}{2 k!}\sum_{l=k}^{n}T_{1}(n,l)(l)_{k}(1+\delta_{l,k})
=\frac{1}{2}\sum_{l=k}^{n}\binom{l}{k}T_{1}(n,l)+\frac{1}{2}T_{1}(n,k).
\end{align*}
Thus we have $T_{1}(n,k;\bold{P})=\frac{1}{2}\sum_{l=k}^{n}\binom{l}{k}T_{1}(n,l)+\frac{1}{2}T_{1}(n,k)$.  By Theorem 5.3, we obtain
\begin{align*}
&\sum_{k=l}^{n}\sum_{m=k}^{n}\sum_{j=l}^{k}\frac{1}{2}\binom{m}{k}\binom{k}{j}T_{1}(n,m)T_{2}(j,l)E_{k-j}\\
&+\sum_{k=l}^{n}\sum_{j=l}^{k}\frac{1}{2}\binom{k}{j}T_{1}(n,k)T_{2}(j,l)E_{k-j}=\delta_{n,l},\\
&\sum_{k=l}^{n}\sum_{j=k}^{n}\sum_{m=l}^{k}\frac{1}{2}\binom{m}{l}\binom{n}{j}T_{2}(j,k)T_{1}(k,m)E_{n-j}\\
&+\sum_{k=l}^{n}\sum_{j=k}^{n}\frac{1}{2}\binom{n}{j}T_{2}(j,k)T_1(k,l)E_{n-j}=\delta_{n,l},\\
& a_{n}=\sum_{k=0}^{n}\sum_{l=k}^{n}\binom{n}{l}T_{2}(l,k)E_{n-l}c_{k} \Longleftrightarrow c_{n}=\sum_{k=0}^{n}\left\{\frac{1}{2}\sum_{l=k}^{n}\binom{l}{k}T_{1}(n,l)+\frac{1}{2}T_{1}(n,k)\right\}a_{k},\\
&a_{n}=\sum_{k=n}^{m}\sum_{l=n}^{k}\binom{k}{l}T_{2}(l,n)E_{k-l}c_{k} \Longleftrightarrow c_{n}=\sum_{k=n}^{m}\left\{\frac{1}{2}\sum_{l=n}^{k}\binom{l}{n}T_{1}(k,l)+\frac{1}{2}T_{1}(k,n)\right\}a_{k}.
\end{align*}

(q) Let $\bold{P}=\left\{(rx+s)_{n}\right\}$, with $ r \ne 0$. Then $(rx+s)_{n} \sim \big(e^{-\frac{s}{r} t}, e^{\frac{t}{r}}-1\big)$, and $T_{2}(n,k;\bold{P})=\sum_{l=k}^{n}\sum_{m=l}^{n}\binom{n}{m}r^{l}(s)_{n-m}S_{1}(m,l)T_{2}(l,k)$ (see [26]). By \eqref{7B}, \eqref{18B} and Theorem 5.1, we see that
\begin{align*}
T_{1}(n,k;\bold{P})&=\frac{1}{k!}\sum_{l=k}^{n}T_{1}(n,l)\big\langle{e^{-\frac{s}{r}t}\big(e^{\frac{t}{r}}-1\big)^{k} \big | x^{l} \big \rangle}\\
&=\frac{1}{k!}\sum_{l=k}^{n}T_{1}(n,l)\Big\langle{(e^{t}-1)^{k} \Big | \Big(\frac{x-s}{r}\Big)^{l} \Big \rangle}\\
&=\sum_{l=k}^{n}T_{1}(n,l)\frac{1}{r^{l}}\frac{1}{k!}\big\langle{(e^{t}-1)^{k} |(x-s)^{l} \big\rangle}\\
&=\sum_{l=k}^{n}T_{1}(n,l)\frac{1}{r^{l}}\sum_{m=k}^{l}S_{2}(m,k)\frac{1}{m!}\big\langle{t^{m} |(x-s)^{l} \big\rangle}\\
&=\sum_{l=k}^{n}\sum_{m=k}^{l}\binom{l}{m}r^{-l}(-s)^{l-m}T_{1}(n,l)S_{2}(m,k).
\end{align*}
Thus we have $T_{1}(n,k;\bold{P})=\sum_{l=k}^{n}\sum_{m=k}^{l}\binom{l}{m}r^{-l}(-s)^{l-m}T_{1}(n,l)S_{2}(m,k)$. From Theorem 5.3, we obtain
\begin{align*}
&a_{n}=\sum_{k=0}^{n}\sum_{l=k}^{n}\sum_{m=l}^{n}\binom{n}{m}r^{l}(s)_{n-m}S_{1}(m,l)T_{2}(l,k)c_{k}\\
& \Longleftrightarrow c_{n}=\sum_{k=0}^{n}\sum_{l=k}^{n}\sum_{m=k}^{l}\binom{l}{m}r^{-l}(-s)^{l-m}T_{1}(n,l)S_{2}(m,k)a_{k},\\
&a_{n}=\sum_{k=n}^{m}\sum_{l=n}^{k}\sum_{m=l}^{k}\binom{k}{m}r^{l}(s)_{k-m}S_{1}(m,l)T_{2}(l,n)c_{k}\\
& \Longleftrightarrow c_{n}=\sum_{k=n}^{m}\sum_{l=n}^{k}\sum_{m=n}^{l}\binom{l}{m}r^{-l}(-s)^{l-m}T_{1}(k,l)S_{2}(m,n)a_{k}.
\end{align*}

(r) Let $\bold{P}=\left\{b_{n}(x)\right\}$ be the sequence of Bernoulli polynomials of the second kind. Then $b_{n}(x) \sim \big(\frac{t}{e^{t}-1}, e^{t}-1\big)$, and $T_{2}(n,k;\bold{P})=\sum_{l=k}^{n}\sum_{m=l}^{n}\binom{n}{m}b_{n-m}S_{1}(m,l)T_{2}(l,k)$,
where $b_{n}=b_{n}(0)$ are the Bernoulli numbers of the second (see [31]). By using \eqref{10B}, \eqref{12B} and Theorem 5.1 and noting $B_{n}(x) \sim (\frac{e^{t}-1}{t},t)$,  we have 
\begin{align*}
T_{1}(n,k;\bold{P})&=\sum_{l=k}^{n}T_{1}(n,l)\Big\langle{\frac{1}{k!}(e^{t}-1)^{k} \Big| \frac{t}{e^{t}-1}x^{l}\Big\rangle}=\sum_{l=k}^{n}T_{1}(n,l)\Big\langle{\frac{1}{k!}(e^{t}-1)^{k} \Big| B_{l}(x)\Big\rangle}\\
&=\sum_{l=k}^{n}T_{1}(n,l)\sum_{m=k}^{l}S_{2}(m,k)\frac{1}{m!}\big\langle{t^{m}|B_{l}(x)\rangle}
=\sum_{l=k}^{n}\sum_{m=k}^{l}\binom{l}{m}B_{l-m}T_{1}(n,l)S_{2}(m,k),
\end{align*}
where $B_{l}$ are the Bernoulli numbers. Hence $T_{1}(n,k;\bold{P})=\sum_{l=k}^{n}\sum_{m=k}^{l}\binom{l}{m}B_{l-m}T_{1}(n,l)S_{2}(m,k)$.
Now, by Theorem 5.3, we get
\begin{align*}
&a_{n}=\sum_{k=0}^{n}\sum_{l=k}^{n}\sum_{m=l}^{n}\binom{n}{m}b_{n-m}S_{1}(m,l)T_{2}(l,k)c_{k}\\
&\Longleftrightarrow c_{n}=\sum_{k=0}^{n}\sum_{l=k}^{n}\sum_{m=k}^{l}\binom{l}{m}B_{l-m}T_{1}(n,l)S_{2}(m,k)a_{k},\\
&a_{n}=\sum_{k=n}^{m}\sum_{l=n}^{k}\sum_{m=l}^{k}\binom{k}{m}b_{k-m}S_{1}(m,l)T_{2}(l,n)c_{k} \\
&\Longleftrightarrow c_{n}=\sum_{k=n}^{m}\sum_{l=n}^{k}\sum_{m=n}^{l}\binom{l}{m}B_{l-m}T_{1}(k,l)S_{2}(m,n)a_{k}.
\end{align*}

(s) Let $\bold{P}=\left\{C_{n}(x;a)\right\}$ be the sequence of Poisson-Charlier polynomials. Then $C_{n}(x;a) \sim \big(e^{a(e^{t}-1)}, a(e^{t}-1)\big)$, with $a \ne 0$, and $T_{2}(n,k;\bold{P})=\sum_{l=k}^{n}\sum_{m=l}^{n}\binom{n}{m}(-1)^{n-m}a^{-m}S_{1}(m,l)T_{2}(l,k)$ (see [31]). By Theorem 5.1, we get
\begin{align*}
T_{1}(n,k;\bold{P})&=\frac{1}{k!}\sum_{l=k}^{n}T_{1}(n,l)a^{k}\langle{(e^{t}-1)^{k}|e^{a(e^{t}-1)}x^{l}\rangle}\\
&=\sum_{l=k}^{n}T_{1}(n,l)a^{k}\sum_{m=k}^{l}\binom{l}{m}\mathrm{Bel}_{l-m}(a)\big\langle{\frac{1}{k!}(e^{t}-1)^{k}|x^{m}\big\rangle}\\
&=\sum_{l=k}^{n}T_{1}(n,l)a^{k}\sum_{m=k}^{l}\binom{l}{m}\mathrm{Bel}_{l-m}(a)\sum_{j=k}^{m}
S_{2}(j,k)\frac{1}{j!}\langle{t^{j}|x^{m}\rangle}\\
&=\sum_{l=k}^{n}\sum_{m=k}^{l}a^{k}\binom{l}{m}\mathrm{Bel}_{l-m}(a)T_{1}(n,l)S_{2}(m,k).
\end{align*}
where $\mathrm{Bel}_{n}(x)$ are the Bell polynomials given by $e^{x(e^{t}-1)}=\sum_{n=0}^{\infty}\mathrm{Bel}_{n}(x)\frac{t^{n}}{n!}$.\\
Thus we have $T_{1}(n,k;\bold{P})=\sum_{l=k}^{n}\sum_{m=0}^{l}a^{k}\binom{l}{m}\mathrm{Bel}_{l-m}(a)T_{1}(n,l)S_{2}(m,k)$.

Now, by Theorem 5.3, we obtain
\begin{align*}
&a_{n}=\sum_{k=0}^{n}\sum_{l=k}^{n}\sum_{m=l}^{n}\binom{n}{m}(-1)^{n-m}a^{-m}S_{1}(m,l)T_{2}(l,k)c_{k} \\
&\Longleftrightarrow c_{n}=\sum_{k=0}^{n}\sum_{l=k}^{n}\sum_{m=0}^{l}a^{k}\binom{l}{m}\mathrm{Bel}_{l-m}(a)T_{1}(n,l)S_{2}(m,k)a_{k},\\
&a_{n}=\sum_{k=n}^{m}\sum_{l=n}^{k}\sum_{m=l}^{k}\binom{k}{m}(-1)^{k-m}a^{-m}S_{1}(m,l)T_{2}(l,n)c _{k}\\
& \Longleftrightarrow c_{n}=\sum_{k=n}^{m}\sum_{l=n}^{k}\sum_{m=0}^{l}a^{n}\binom{l}{m}\mathrm{Bel}_{l-m}(a)T_{1}(k,l)S_{2}(m,n)a_{k}.
\end{align*}

(t) Let $\bold{P}=\left\{p_{n}(x)\right\}$, with $p_{n}(x)=\sum_{k=0}^{n}B_{k}(x)B_{n-k}(x)$. This is not a Sheffer sequence. Here we would like to describe how to determine $T_{1}(n,k;\bold{P})$.
For this, we recall from [17] that
\begin{align}
p_{n}(x)=\frac{2}{n+2}\sum_{m=0}^{n-2}\binom{n+2}{m}B_{n-m}B_m(x)+(n+1)B_n(x).\label{2F}
\end{align}
On the one hand, by (o) we have
\begin{align}
x^{[n]}=\sum_{m=0}^{n}\gamma_{m}B_{m}(x), \label{3F}
\end{align}
where $\gamma_{m}=\sum_{l=m}^{n}\frac{1}{l+1}\binom{l+1}{m}T_{1}(n,l)$.\\
On the other hand, by definition and from \eqref{2F} we also have
\begin{align}
x^{[n]}&=\sum_{k=0}^{n}T_{1}(n,k;\bold{P})p_{k}(x)\nonumber\\
&=\sum_{k=2}^{n}T_{1}(n,k;\bold{P})\frac{2}{k+2}\sum_{m=0}^{k-2}\binom{k+2}{m}B_{k-m}B_m(x) \nonumber\\
&\quad\quad+\sum_{k=0}^{n}T_{1}(n,k;\bold{P})(k+1)B_k(x)\label{4F}\\
&=\sum_{m=0}^{n-2}\sum_{k=m+2}^{n}\epsilon_{m,k}T_{1}(n,k;\bold{P})B_{m}(x)\nonumber\\
&\quad\quad+\sum_{m=0}^{n}(m+1)T_{1}(n,m;\bold{P})B_m(x),\nonumber
\end{align}
where $\epsilon_{m,k}=\frac{2}{k+2}\binom{k+2}{m}B_{k-m},\quad (0 \le m \le n-2, m+2 \le k \le n)$.\\
By comparing \eqref{3F} and \eqref{4F}, we have
\begin{align}
&\gamma_{m}=(m+1)T_{1}(n,m;\bold{P})+\sum_{k=m+2}^{n}\epsilon_{m,k}T_{1}(n,k;\bold{P}),
\quad(0 \le m \le n), \label{5F}
\end{align}
where we understand that the sum is zero for $m=n-1$ and $m=n$. Let $\Gamma$ be the column vector consisting of $\gamma_{0}, \gamma_{1}, \dots, \gamma_{n}$, and let $S$ be the column vector consisting of $T_{1}(n,0;\bold{P}), T_{1}(n,1;\bold{P}), \dots, T_{1}(n,n;\bold{P})$. Then,  in matrix form,  \eqref{5F} can be written as  $\Gamma=AS$, where $A$ is an $(n+1) \times (n+1)$ upper triangular matrix with diagonal entries $1,2, \dots , n+1$. Thus $S=A^{-1}\Gamma$.

\section{Conclusion}
In this paper, we studied central factorial numbers of both kinds associated with sequences of polynomilas by using umbral calculus techniques, which was motivated by the observation that many important special numbers appear in the expansions of some polynomials in terms of central factorials and vice versa. Our results were illustrated by twenty examples both for the second kind and for the first kind. In all cases of these examples, some interesting orthogonality and  inverse relations were obtained as immediate consequences of our results. The novelty of this paper is that it is the first paper which studies the central factorial numbers of both kinds associated with any sequence of polynomials in a unified and systematic way with the help of umbral calculus. \\
It is one of our future projects to study special numbers and polynomials by making use of various tools, including combinatorial methods, generating functions, umbral calculus, $p$-adic analysis, probability theory, differential equations, analytic number theory, operator theory, special functions and so on.


\end{document}